\documentclass[10pt,twoside,reqno]{amsart}

\usepackage[margin=1in]{geometry}

\usepackage{mathtools}
\usepackage{enumitem}
\usepackage{mathrsfs}
\usepackage{dsfont}
\usepackage{url}
\usepackage{booktabs}

\usepackage[math]{cellspace}
\newtheorem{theorem}{\sc\hspace{8pt} Theorem}
\newtheorem{corollary}{\sc\hspace{8pt} Corollary}
\newtheorem*{acknowledgments}{Acknowledgments}

\newcommand{\abs}[1]{\lvert#1\rvert}

\DeclareMathOperator*{\im}{Im}
\DeclareMathOperator*{\re}{Re}

\begin{document}

\title[The level crossings of random sums]{The level crossings of random sums}

\author[Christopher Corley]{Christopher Corley}

\address{Department of Mathematics, The University of Tennessee at Chattanooga, 304 Lupton Hall (Dept. 6956), 615 McCallie Avenue, Chattanooga, Tennessee 37403, United States of America}

\email{Christopher-Corley@utc.edu}

\author[Andrew Ledoan]{Andrew Ledoan}

\address{Department of Mathematics, The University of Tennessee at Chattanooga, 304 Lupton Hall (Dept. 6956), 615 McCallie Avenue, Chattanooga, Tennessee 37403, United States of America}

\email{Andrew-Ledoan@utc.edu}

\subjclass[2010]{Primary 30C15; Secondary 30B20, 60B99, 60E05, 60G50, 62H10}

\keywords{Density function; Level crossing; Multivariate normal distribution; Rice formula; Random field; Random summation; Standard Brownian motion}

\begin{abstract}
Let $\{\eta_{j}\}_{j = 0}^{N}$ be a sequence of independent and identically distributed complex normal random variables with mean zero and variances $\{\sigma_{j}^{2}\}_{j = 0}^{N}$. Let $\{f_{j} (z)\}_{j = 0}^{N}$ be a sequence of holomorphic functions that are real-valued on the real line. The purpose of the present study is that of examining the number of times that the random sum $\sum_{j = 0}^{N} \eta_{j} f_{j} (z)$ crosses the complex level $\boldsymbol{K} = K_{1} + i K_{2}$, where $K_{1}$ and $K_{2}$ are constants independent of $z$. More specifically, we establish an exact formula for the expected density function for the complex zeros. We then reformulate the problem in terms of successive observations of a Brownian motion. We further answer the basic question about the expected number of complex zeros for coefficients of nonvanishing mean values.
\end{abstract}

\maketitle

\thispagestyle{empty}

\section{Introduction and statement of results}

An exact formula for the expected number of real zeros of a random polynomial was obtained by Kac \cite{Kac1943} under independent and identically distributed (i.i.d.) real standard normal coefficients. For complex coefficients, Dunnage \cite{Dunnage1968, Dunnage1970} gave some estimates for the number of real zeros. For complex zeros, the expected density of zeros was studied by Shepp and Vanderbei \cite{SheppVanderbei1995} for i.i.d. real standard normal coefficients and generalized by Ibragimov and Zeitouni \cite{IbragimovZeitouni1997} for a wider class of distributions of coefficients. Relevant to these investigations is the work of Kostlan \cite{Kostlan1993}. The expected density was dealt with, also, by Hammersley \cite{Hammersley1956}, Edelman and Kostlan \cite{EdelmanKostlan1995}, and Farahmand and Grigorash \cite{FarahmandGrigorash1999}. Vanderbei \cite{Vanderbei2015} generalized the work in \cite{SheppVanderbei1995} to random sums with holomorphic functions that are real-valued on the real line as basis functions. Motivated by the studies conducted by Vanderbei \cite{Vanderbei2015} and Farahmand \cite{Farahmand1997}, the present authors \cite{CorleyLedoan2020} studied the number of times that these random sums cross a complex level and obtained certain results on the level crossings. The aim of the present paper is to extend these results in the following direction.

Let $z$ be the complex variable $x + i y$. Let $\{a_{j}\}_{j = 0}^{N}$ and $\{b_{j}\}_{j = 0}^{N}$ be sequences of mutually i.i.d. real normal random variables defined on the complete probability space $(\Omega, \mathscr{F}, \operatorname{Prob})$ with mean zero and variances $\{\sigma_{a_{j}}^{2}\}_{j = 0}^{N}$ and $\{\sigma_{b_{j}}^{2}\}_{j = 0}^{N}$, so that $\mathscr{F}$ is a $\sigma$-field of subsets of $\Omega$ and $\operatorname{Prob}$ is a probability measure on $(\Omega, \mathscr{F})$. Assume all sub $\sigma$-fields contain all sets of measure zero. Then, let $\{\eta_{j}\}_{j = 0}^{N}$ be a sequence of i.i.d. complex normal random variables with density $e^{- z \bar{z}} / \pi$ and given by $\eta_{j} = a_{j} + i b_{j}$ for $0 \leq j \leq N$. Let, further, $\{f_{j} (z)\}_{j = 0}^{N}$ be a sequence of holomorphic functions $f_{j} (z) = u_{j} (x, y) + i v_{j} (x, y)$ for $0 \leq j \leq N$ that are real-valued on the real line. By the Schwarz reflection principle, $\overline{f_{j} (z)} = f_{j} (\overline{z})$ for $0 \leq j \leq N$. Then, define
\begin{equation} \label{eq-1}
S_{N} (z)
 = \eta_{0} f_{0} (z) + \eta_{1} f_{1} (z) + \dots + \eta_{N} f_{N} (z)
 = \sum_{j = 0}^{N} \eta_{j} f_{j} (z).
\end{equation}

If, for each compact subset $T$ of the complex plane, $N_{\boldsymbol{K}}^{S_{N}} (T)$ denotes the random number of complex zeros, counted with multiplicity, in $T$ of $S_{N} (z)$ that cross the complex level $\boldsymbol{K} = K_{1} + i K_{2}$, where $K_{1}$ and $K_{2}$ are constants independent of $z$, then from \cite{CorleyLedoan2020}, with probability one, the expected density $h_{N, \boldsymbol{K}} (z)$ of the complex zeros of
\begin{equation} \label{eq-2}
S_{N} (z)
 = \boldsymbol{K}
\end{equation}
is given by
\begin{equation} \label{eq-3}
E (N_{\boldsymbol{K}}^{S} (T))
 = \int_{T} h_{N, \boldsymbol{K}} (z) \, dz.
\end{equation}
The explicit derivation of $h_{N, \boldsymbol{K}} (z)$ constitutes the primary reason for studying the complex zeros of \eqref{eq-2}. The main device for treating $h_{N, \boldsymbol{K}} (z)$ throughout the complex plane is the Rice formula. This remarkable result provides a representation for the expected number of zeros of certain random fields. It is reproduced below from \cite[Theorem 6.2, pp. 163-164]{AzaisWschebor2009}. (See, also, \cite[Theorem 11.2.3, Corollary 11.2.4, pp. 269-271]{AdlerTaylor2007} and \cite[Theorem 2.1, p. 256]{AzaisWschebor2005}.)
\begin{theorem} \label{th-1}
Let $Z \colon U \rightarrow \mathds{R}^{N}$ be a random field, let $U$ be an open subset of $\mathds{R}^{N}$, and let $u \in \mathds{R}^{N}$ be a fixed point in the codomain. Assume the following conditions are satisfied with probability one:
\begin{enumerate}[label=(\roman{*})]
\item $Z$ is normal.
\item Almost surely the function $t \rightsquigarrow Z (t)$ is of class $C^{1}$.
\item For each $t \in U$, $Z (t)$ has a nondegenerate distribution---i.e., $\operatorname{Var} (Z (t)) \succ 0$.
\item For each $u \in \mathds{R}^{N}$, $\operatorname{Prob} (\exists \, t \in U \colon Z (t) = u, \det (Z' (t)) = 0) = 0$.
\end{enumerate}
If $N_{u}^{Z} (B)$ denotes the number of zeros of $Z (t) = u$ that belong to the Borel subset $B \subset U$, then one has
\begin{equation} \label{eq-4}
E (N_{u}^{Z} (B))
 = \int_{B} E (\abs{\det (Z' (t))} \mid Z (t) = u) \, p_{Z (t)} (u) \, dt,
\end{equation}
where $p_{Z (t)} (u)$ is the probability density function of $Z (t)$ at $u$. If $B$ is compact, then both sides of \eqref{eq-4} are finite.
\end{theorem}

The function $Z$ in \eqref{eq-4} is defined on $\mathds{R}^{N}$. In our application, we need to find the real and complex zeros of \eqref{eq-2}---i.e., the real zeros of $\re (S_{N} (x + i y)) = K_{1}$ and $\im (S_{N} (x + i y)) = K_{2}$ for $(x, y) \in \mathds{R}^{2}$. The conditions $(i)$--$(iv)$ are easy to check. Formula \eqref{eq-4} is interesting. It shows that $h_{N, \boldsymbol{K}} (z)$, as defined by \eqref{eq-3}, can be expressed through a conditioned mean function of a quadratic form of i.i.d. real normal random variables conditioned on certain linear combinations.

\begin{theorem} \label{th-2}
Provided all the conditions imposed on $S_{N} (z)$ in \eqref{eq-1} and $T$ are satisfied, then for all $N > 1$ one has
\begin{equation*}
\begin{split}
& h_{N, \boldsymbol{K}} (z)
 = \frac{1}{2 \pi D_{0} (z)} \exp \left(-\frac{K_{1}^{2} Y_{3} (z) + K_{2}^{2} Y_{1} (z) - 2 K_{1} K_{2} Y_{2} (z)}{2 D_{0} (z)^{2}}\right) \\ & \hspace{0.5em} \times \left\{D_{3} (z) - \frac{\abs{D_{1} (z)}^{2}}{D_{0} (z)} \left(\frac{Y_{2} (z) + Y_{3} (z)}{D_{0} (z)} - \frac{(K_{1} Y_{3} (z) - K_{2} Y_{2} (z)) (K_{1} (Y_{2} (z) + Y_{3} (z)) - K_{2} (Y_{1} (z) + Y_{2} (z)))}{D_{0} (z)^{3}}\right)\right. \\ & \hspace{3.5em} - \frac{\abs{D_{2} (z)}^{2}}{D_{0} (z)} \left(\frac{Y_{1} (z) + Y_{2} (z)}{D_{0} (z)} - \frac{(K_{1} Y_{2} (z) - K_{2} Y_{1} (z)) (K_{1} (Y_{2} (z) + Y_{3} (z)) - K_{2} (Y_{1} (z) + Y_{2} (z)))}{D_{0} (z)^{3}}\right) \\ & \hspace{6.5em} \left. + \frac{\abs{D_{1} (z) + i D_{2} (z)}^{2}}{D_{0} (z)} \left(\frac{Y_{2} (z)}{D_{0} (z)} - \frac{(K_{1} Y_{3} (z) - K_{2} Y_{2} (z)) (K_{1} Y_{2} (z) - K_{2} Y_{1} (z))}{D_{0} (z)^{3}}\right)\right\},
\end{split}
\end{equation*}
where
\begin{equation*}
\begin{array}{c@{\qquad}c}
\displaystyle Y_{1} (z)
 = \sum_{j = 0}^{N} (\sigma_{a_{j}}^{2} u_{j}^{2} + \sigma_{b_{j}}^{2} v_{j}^{2}), \qquad
\displaystyle Y_{2} (z)
 = \sum_{j = 0}^{N} (\sigma_{a_{j}}^{2} - \sigma_{b_{j}}^{2}) u_{j} v_{j}, \qquad
\displaystyle Y_{3} (z)
 = \sum_{j = 0}^{N} (\sigma_{b_{j}}^{2} u_{j}^{2} + \sigma_{a_{j}}^{2} v_{j}^{2}),
\end{array}
\end{equation*}
and
\begin{equation*}
\begin{array}{c@{\qquad}c}
\displaystyle D_{0} (z)
 = \sqrt{Y_{1} (z) Y_{3} (z) - Y_{2}^{2} (z)}, \qquad
\displaystyle D_{1} (z)
 = \sum_{j = 0}^{N} (\sigma_{a_{j}}^{2} u_{j} - i \sigma_{b_{j}}^{2} v_{j}) \left(\frac{\partial u_{j}}{\partial x} + i \frac{\partial v_{j}}{\partial x}\right), \\ [3ex]
\displaystyle D_{2} (z)
 = \sum_{j = 0}^{N} (\sigma_{b_{j}}^{2} u_{j} - i \sigma_{a_{j}}^{2} v_{j}) \left(\frac{\partial u_{j}}{\partial x} + i \frac{\partial v_{j}}{\partial x}\right), \qquad
\displaystyle D_{3} (z)
 = \sum_{j = 0}^{N} (\sigma_{a_{j}}^{2} + \sigma_{b_{j}}^{2}) \left(\left(\frac{\partial u_{j}}{\partial x}\right)^{2} + \left(\frac{\partial v_{j}}{\partial x}\right)^{2}\right).
\end{array}
\end{equation*}
\end{theorem}

In relation to the work in \cite{CorleyLedoan2020}, observe that when $\sigma_{a_{j}}^{2} = \sigma_{b_{j}}^{2} = \sigma^{2}$ for $0 \leq j \leq N$
\begin{equation*}
\begin{split}
\begin{array}{c@{\qquad}c}
\displaystyle Y_{1} (z)
 = Y_{3} (z)
 = \sigma^{2} B_{0} (z), \qquad
\displaystyle Y_{2} (z)
 = 0,
\end{array}
\end{split}
\end{equation*}
and
\begin{equation*}
\begin{split}
\begin{array}{c@{\qquad}c}
\displaystyle D_{0} (z)
 = \sigma^{2} B_{0} (z), \qquad
\displaystyle D_{1} (z)
 = D_{2} (z)
 = \sigma^{2} B_{1} (z), \qquad
\displaystyle D_{3} (z)
 = 2 \sigma^{2} B_{2} (z).
\end{array}
\end{split}
\end{equation*}
Then
\begin{equation*}
\abs{D_{1} (z) + i D_{2} (z)}^{2}
 = \abs{D_{1} (z)}^{2} + \abs{D_{2} (z)}^{2}
 = 2 \sigma^{4} \abs{B_{1} (z)}^{2}.
\end{equation*}
The following result is obtained by using these substitutions in Theorem \ref{th-2}.

\begin{theorem}
If $\sigma_{a_{j}}^{2} = \sigma_{b_{j}}^{2} = \sigma^{2}$ for $0 \leq j \leq N$, then for all $N > 1$ one has
\begin{equation*}
\begin{split}
h_{N, \boldsymbol{K}} (z)
 = \frac{1}{\pi B_{0} (z)} \exp \left(- \frac{K_{1}^{2} + K_{2}^{2}}{2 \sigma^{2} B_{0} (z)}\right) \left\{B_{2} (z) - \left(\frac{\abs{B_{1} (z)}}{B_{0} (z)}\right)^{2} \left(B_{0} (z) - \frac{K_{1}^{2} + K_{2}^{2}}{2 \sigma^{2}}\right)\right\},
\end{split}
\end{equation*}
where
\begin{equation*}
\begin{split}
\begin{array}{c@{\qquad}c}
\displaystyle B_{0} (z)
 = \sum_{j = 0}^{N} \abs{f_{j} (z)}^{2}, \qquad
\displaystyle B_{1} (z)
 = \sum_{j = 0}^{N} \overline{f_{j} (z)} f_{j}^{\prime} (z), \qquad
\displaystyle B_{2} (z)
 = \sum_{j = 0}^{N} \abs{f_{j}^{\prime} (z)}^{2}.
\end{array}
\end{split}
\end{equation*}
\end{theorem}

Consequently, when $\sigma^{2}$ is set to be one, Theorem 1 in \cite{CorleyLedoan2020} is recovered. Further, if $\boldsymbol{K}$ is the zero vector, Corollary 3 in \cite{CorleyLedoan2020} is recovered. This latter result was proved independently by Yeager \cite{Yeager2016} and one of the authors \cite{Ledoan2016}.

In addition, the following result follows from Theorem \ref{th-2}.

\begin{corollary} \label{cor-1}
If the vector $\boldsymbol{K}$ is restricted to a circle of radius $K > 0$, then for all $N > 1$ one has
\begin{equation*}
\begin{split}
h_{N, \boldsymbol{K}} (z)
 & = \frac{1}{2 \pi D_{0} (z)} \exp \left(-\frac{K^{2} (Y_{1} (z) - Y_{2} (z) + Y_{3} (z))}{2 D_{0} (z)^{2}}\right) \\ & \hspace{3em} \times \left\{D_{3} (z) - \frac{\abs{D_{1} (z)}^{2}}{D_{0} (z)} \left(\frac{Y_{2} (z) + Y_{3} (z)}{D_{0} (z)} - \frac{K^{2} (Y_{2} (z) - Y_{3} (z)) (Y_{1} (z) - Y_{2} (z) - 1)}{D_{0} (z)^{3}}\right)\right. \\ & \hspace{6em} - \frac{\abs{D_{2} (z)}^{2}}{D_{0} (z)} \left(\frac{Y_{1} (z) + Y_{2} (z)}{D_{0} (z)} - \frac{K^{2} (Y_{1} (z) - Y_{2} (z)) (Y_{2} (z) - Y_{3} (z) + 1)}{D_{0} (z)^{3}}\right) \\ & \hspace{9em} \left. + \frac{\abs{D_{1} (z) + i D_{2} (z)}^{2}}{D_{0} (z)} \left(\frac{Y_{2} (z)}{D_{0} (z)} - \frac{K^{2} (Y_{1} (z) - Y_{2} (z)) (Y_{2} (z) - Y_{3} (z))}{D_{0} (z)^{3}}\right)\right\}.
\end{split}
\end{equation*}
\end{corollary}

A special case of Corollary \ref{cor-1} follows.

\begin{corollary}
If $\boldsymbol{K}$ is the zero vector, then for all $N > 1$ one has
\begin{equation*}
\begin{split}
h_{N, \boldsymbol{0}} (z)
 & = \frac{D_{0} (z)^{2} D_{3} (z) - \abs{D_{1} (z)}^{2} (Y_{2} (z) + Y_{3} (z)) - \abs{D_{2} (z)}^{2} (Y_{1} (z) + Y_{2} (z)) + \abs{D_{1} (z) + i D_{2} (z)}^{2} Y_{2} (z)}{2 \pi D_{0} (z)^{3}}.
\end{split}
\end{equation*}
\end{corollary}

The proof of Theorem \ref{th-2}, in the spirit of the method credited to Ibragimov and Zeitouni \cite{IbragimovZeitouni1997}, is presented in Section \ref{sec-2}. In relation to the works of Rezakhah and Shemehsavar \cite{RezakhahShemehsavar2005} and Rezakhah and Soltani \cite{RezakhahSoltani2003}, an application of Theorem \ref{th-2} entailing a sequence of successive observations of a Brownian motion is given in Section \ref{sec-3}. Theorem \ref{th-2} can be extended further by letting the $a_{j}$ and $b_{j}$ to be mutually i.i.d. real random variables such that $a_{j} \sim \mathcal{N} (\mu_{a_{j}}, \sigma_{a_{j}}^{2})$ and $b_{j} \sim \mathcal{N} (\mu_{b_{j}}, \sigma_{b_{j}}^{2})$ for $0 \leq j \leq N$.

Toward this end, let
\begin{equation*}
\begin{split}
Y_{1}^{\ast} (z)
 &= \sum_{j = 0}^{N} (\sigma_{a_{j}}^{2} u_{j}^{2} + \sigma_{b_{j}}^{2} v_{j}^{2}) - \left(\sum_{j = 0}^{N} (\mu_{a_{j}} u_{j} - \mu_{b_{j}} v_{j})\right)^{2}, \\
Y_{2}^{\ast} (z)
 & = \sum_{j = 0}^{N} (\sigma_{a_{j}}^{2} - \sigma_{b_{j}}^{2}) u_{j} v_{j} - \left(\sum_{j = 0}^{N} (\mu_{a_{j}} u_{j} - \mu_{b_{j}} v_{j})\right) \left(\sum_{j = 0}^{N} (\mu_{a_{j}} v_{j} + \mu_{b_{j}} u_{j})\right), \\
Y_{3}^{\ast} (z)
 & = \sum_{j = 0}^{N} (\sigma_{a_{j}}^{2} v_{j}^{2} + \sigma_{b_{j}}^{2} u_{j}^{2}) - \left(\sum_{j = 0}^{N} (\mu_{a_{j}} v_{j} + \mu_{b_{j}} u_{j})\right)^{2}.
\end{split}
\end{equation*}
Then define
\begin{equation*}
M (z)
 = \sum_{j = 0}^{N} E (\eta_{j}) f_{j}^{\prime} (z)
\end{equation*}
and
\begin{equation*}
\begin{split}
D_{0}^{\ast} (z)
 & = \sqrt{Y_{1}^{\ast} (z) Y_{3}^{\ast} (z) - Y_{2}^{\ast} (z)^{2}}, \\
D_{2}^{\ast} (z)
 & = \sum_{j = 0}^{N} (B_{j, 2} (z) - i A_{j, 2} (z)) f_{j}^{\prime} (z),
\end{split}
\qquad
\begin{split}
D_{1}^{\ast} (z)
 & = \sum_{j = 0}^{N} (A_{j, 1} (z) - i B_{j, 1} (z)) f_{j}^{\prime} (z), \\
D_{3}^{\ast} (z)
 & = \sum_{j = 0}^{N} (\sigma_{a_{j}}^{2} + \sigma_{b_{j}}^{2}) \abs{f_{j}^{\prime} (z)}^{2},
\end{split}
\end{equation*}
where
\begin{equation*}
\begin{split}
A_{j, 1} (z)
 & = \sigma_{a_{j}}^{2} u_{j} - \mu_{a_{j}} E (X_{1}), \\
B_{j, 1} (z)
 & = \sigma_{b_{j}}^{2} v_{j} + \mu_{b_{j}} E (X_{1}),
\end{split}
\qquad
\begin{split}
A_{j, 2} (z)
 & = \sigma_{a_{j}}^{2} v_{j} - \mu_{a_{j}} E (X_{2}), \\
B_{j, 2} (z)
 & = \sigma_{b_{j}}^{2} u_{j} - \mu_{b_{j}} E (X_{2}).
\end{split}
\end{equation*}
The following theorem is proved in Section \ref{sec-4}.

\begin{theorem} \label{th-4}
Provided all the conditions imposed on $S_{N} (z)$ in \eqref{eq-1} and $T$ are satisfied, then for all $N > 1$ one has
\begin{equation*}
\begin{split}
& h_{N, \boldsymbol{K}} (z)
 = \frac{1}{2 \pi D_{0}^{\ast} (z)} \\ & \hspace{0.5em} \times \exp \left(-\frac{(K_{1} - E (X_{1}))^{2} Y_{3}^{\ast} (z) + (K_{2} - E (X_{2}))^{2} Y_{1}^{\ast} (z) - 2 (K_{1} - E (X_{1})) (K_{2} - E (X_{2})) Y_{2}^{\ast} (z)}{2 D_{0}^{\ast} (z)^{2}}\right) \\ & \hspace{0.5em} \times \left\{D_{3}^{\ast} (z) - \frac{\abs{D_{1}^{\ast} (z)}^{2}}{D_{0}^{\ast} (z)^{2}} \left(Y_{3}^{\ast} (z) - \frac{((K_{1} - E (X_{1})) Y_{3}^{\ast} (z) - (K_{2} - E (X_{2})) Y_{2}^{\ast} (z))^{2}}{D_{0}^{\ast} (z)^{2}}\right) \right. \\ & \hspace{0.5em} \left. - \frac{\abs{D_{2}^{\ast} (z)}^{2}}{D_{0}^{\ast} (z)^{2}} \left(Y_{1}^{\ast} (z) - \frac{((K_{1} - E (X_{1})) Y_{2}^{\ast} (z) - (K_{2} - E (X_{2})) Y_{1}^{\ast} (z))^{2}}{D_{0}^{\ast} (z)^{2}}\right) \right. \\ & \hspace{0.5em} \left. + \left(\frac{\abs{D_{1}^{\ast} (z) + i D_{2}^{\ast} (z)}^{2} - \abs{D_{1}^{\ast} (z)}^{2} - \abs{D_{2}^{\ast} (z)}^{2}}{D_{0}^{\ast} (z)^{2}}\right) \right. \\ &\hspace{0.5em} \left. \times \left(Y_{2}^{\ast} (z) - \frac{((K_{1} - E (X_{1})) Y_{3}^{\ast} (z) - (K_{2} - E (X_{2})) Y_{2}^{\ast} (z)) ((K_{1} - E (X_{1})) Y_{2}^{\ast} (z) - (K_{2} - E (X_{2})) Y_{1}^{\ast} (z))}{D_{0}^{\ast} (z)^{2}}\right) \right. \\ & \hspace{0.5em} \left. - \left(\frac{\abs{M (z) + D_{1}^{\ast} (z)}^{2} - \abs{M (z)}^{2} - \abs{D_{1}^{\ast} (z)}^{2}}{D_{0}^{\ast} (z)^{2}}\right) ((K_{1} - E (X_{1})) Y_{3}^{\ast} (z) - (K_{2} - E (X_{2})) Y_{2}^{\ast} (z)) \right. \\ & \hspace{0.5em} \left. + \left(\frac{\abs{M (z) + i D_{2}^{\ast} (z)}^{2} - \abs{M (z)}^{2} - \abs{D_{2}^{\ast} (z)}^{2}}{D_{0}^{\ast} (z)^{2}}\right) ((K_{1} - E (X_{1})) Y_{2}^{\ast} (z) - (K_{2} - E (X_{2})) Y_{1}^{\ast} (z))\right\}.
\end{split}
\end{equation*}
\end{theorem}

Several consequences of Theorem \ref{th-4} are of special interest. These are derived in Section \ref{sec-5}.

\section{Proof of Theorem \ref{th-2}} \label{sec-2}

The proof of Theorem \ref{th-2} starts with the decomposition
\begin{equation*}
S_{N} (z)
 = X_{1} + i X_{2},
\end{equation*}
where
\begin{equation*}
X_{1}
 = \sum_{j = 0}^{N} (a_{j} u_{j} - b_{j} v_{j}), \qquad
X_{2}
 = \sum_{j = 0}^{N} (a_{j} v_{j} + b_{j} u_{j}).
\end{equation*}
If the column vector
\begin{equation*}
\boldsymbol{X}
 = (X_{1}, X_{2})^{\prime}
\end{equation*}
genuinely represents a two-dimensional random field, then for $z = x + i y$ the Jacobian matrix of the random transformation $(x, y) \rightarrow (X_{1}, X_{2})$ is
\begin{equation*}
\nabla \boldsymbol{X}
 =
\begin{pmatrix*}[l]
\displaystyle \sum_{j = 0}^{N} \left(a_{j} \frac{\partial u_{j}}{\partial x} - b_{j} \frac{\partial v_{j}}{\partial x}\right) & \displaystyle \sum_{j = 0}^{N} \left(a_{j} \frac{\partial v_{j}}{\partial x} + b_{j} \frac{\partial u_{j}}{\partial x}\right) \\ \addlinespace[1mm]
\displaystyle \sum_{j = 0}^{N} \left(-a_{j} \frac{\partial v_{j}}{\partial x} - b_{j} \frac{\partial u_{j}}{\partial x}\right) & \displaystyle \sum_{j = 0}^{N} \left(a_{j} \frac{\partial u_{j}}{\partial x} - b_{j} \frac{\partial v_{j}}{\partial x}\right)
\end{pmatrix*}
\end{equation*}
and
\begin{align} \label{eq-5}
\det (\nabla \boldsymbol{X})
 & = \sum_{j = 0}^{N} \sum_{k = 0}^{N} \left((a_{j} a_{k} + b_{j} b_{k}) \left(\frac{\partial u_{j}}{\partial x} \frac{\partial u_{k}}{\partial x} + \frac{\partial v_{j}}{\partial x} \frac{\partial v_{k}}{\partial x}\right) + (a_{j} b_{k} - b_{j} a_{k}) \left(\frac{\partial v_{j}}{\partial x} \frac{\partial u_{k}}{\partial x} - \frac{\partial u_{j}}{\partial x} \frac{\partial v_{k}}{\partial x}\right)\right) \nonumber \\ \nonumber
 & = \sum_{j = 0}^{N} (a_{j}^{2} + b_{j}^{2}) \left(\left(\frac{\partial u_{j}}{\partial x}\right)^{2} + \left(\frac{\partial v_{j}}{\partial x}\right)^{2}\right) + \sum_{j = 0}^{N} \sum_{\substack{k = 0 \\ k \neq j}}^{N} \left((a_{j} a_{k} + b_{j} b_{k}) \left(\frac{\partial u_{j}}{\partial x} \frac{\partial u_{k}}{\partial x} + \frac{\partial v_{j}}{\partial x} \frac{\partial v_{k}}{\partial x}\right)\right. \nonumber \\ & \hspace{22em} \left. + (a_{j} b_{k} - b_{j} a_{k}) \left(\frac{\partial v_{j}}{\partial x} \frac{\partial u_{k}}{\partial x} - \frac{\partial u_{j}}{\partial x} \frac{\partial v_{k}}{\partial x}\right)\right).
\end{align}
The evaluation of $h_{N, \boldsymbol{K}} (z)$ leads to the computation of the expected value of a quadratic form $\det (\nabla \boldsymbol{X})$ of i.i.d. real normal random variables conditioned on two linear combinations. Observe that $\det (\nabla \boldsymbol{X})$ is always nonnegative. Since $N$ is fixed, $T$ contains not more than a finite number of zeros of
\begin{equation} \label{eq-6}
\boldsymbol{X}
 = \boldsymbol{K},
\end{equation}
where
\begin{equation} \label{eq-7}
\boldsymbol{K}
 = (K_{1}, K_{2})^{\prime}.
\end{equation}
Since the set of zeros of \eqref{eq-6} is of measure zero, assume the boundary $\partial T$ does not contain any zeros of \eqref{eq-6} and $T$ does not contain any such zeros such that $\det (\nabla \boldsymbol{X}) = 0$. Theorem \ref{th-1} applies. Thus,
\begin{equation} \label{eq-8}
h_{N, \boldsymbol{K}} (z)
 = E (\det (\nabla \boldsymbol{X}) \mid \boldsymbol{X} = \boldsymbol{K}) \, p_{X_{1}, X_{2}} (\boldsymbol{K}^{\prime}),
\end{equation}
where $ p_{X_{1}, X_{2}} (\boldsymbol{K}^{\prime})$ denotes the probability density of the random vector $\boldsymbol{X}$ at $\boldsymbol{K}^{\prime}$. By \eqref{eq-6}, and since $X_{1}$ and $X_{2}$ are linear forms with respect to $\eta_{j}$ for $0 \leq j \leq N$, $h_{N, \boldsymbol{K}} (z)$ is the conditional mean of a quadratic form with respect to $\eta_{j}$ for $0 \leq j \leq N$. This form can be calculated in terms of components by means of multivariate analysis.

Based on the assumption that the scalar random variables are independent and normally distributed, the multivariate random vectors
\begin{equation*}
\boldsymbol{a}
 = (a_{0}, a_{1}, \ldots, a_{N})^{\prime}, \qquad
\boldsymbol{b}
 = (b_{0}, b_{1}, \ldots, b_{N})^{\prime}
\end{equation*}
are such that
\begin{equation} \label{eq-9}
\operatorname{Cov} (\boldsymbol{a}, \boldsymbol{b} \mid \boldsymbol{X} = \boldsymbol{K})
 = 
\begin{pmatrix*}[l]
\boldsymbol{\Sigma}_{\boldsymbol{a} \boldsymbol{a}, \boldsymbol{X}} & \boldsymbol{\Sigma}_{\boldsymbol{a} \boldsymbol{b}, \boldsymbol{X}} \\ \addlinespace[1mm]
\boldsymbol{\Sigma}_{\boldsymbol{b} \boldsymbol{a}, \boldsymbol{X}} & \boldsymbol{\Sigma}_{\boldsymbol{b} \boldsymbol{b}, \boldsymbol{X}}
\end{pmatrix*}.
\end{equation}
The elements of this covariance matrix are computed using
\begin{equation} \label{eq-10}
\boldsymbol{\Sigma}_{\boldsymbol{a} \boldsymbol{b}, \boldsymbol{X}}
 = \boldsymbol{\Sigma}_{\boldsymbol{a} \boldsymbol{b}} - \boldsymbol{\Sigma}_{\boldsymbol{a} \boldsymbol{X}} \boldsymbol{\Sigma}_{\boldsymbol{X} \boldsymbol{X}}^{-1} \boldsymbol{\Sigma}_{\boldsymbol{X} \boldsymbol{b}}
\end{equation}
and the corresponding expression
\begin{equation*}
\boldsymbol{\Sigma}_{\boldsymbol{a} \boldsymbol{b}}
 = E ((\boldsymbol{a}- E (\boldsymbol{a})) (\boldsymbol{b} - E (\boldsymbol{b}))^{\prime}).
\end{equation*}
Since the distribution of $a_{j}$ and $b_{j}$ is central for $0 \leq j \leq N$, $E (\boldsymbol{a}) = 0$ and $E (\boldsymbol{b}) = 0$. Clearly,
\begin{equation} \label{eq-11}
\boldsymbol{\Sigma}_{\boldsymbol{a} \boldsymbol{b}}
 = E (\boldsymbol{a} \boldsymbol{b}^{\prime}).
\end{equation}
Thusly, the conditional expected values are expressed in terms of unconditional expected values and covariances.

Then, if $E (X_{1}) = 0$ and $E (X_{2}) = 0$, $E (\boldsymbol{X}) = 0$, whence, by \eqref{eq-11},
\begin{equation} \label{eq-12}
\boldsymbol{\Sigma}_{\boldsymbol{X} \boldsymbol{X}}
 =
\begin{pmatrix*}[l]
E (X_{1} X_{1})	& E (X_{1} X_{2}) \\ \addlinespace[1mm]
E (X_{2} X_{1})	& E (X_{2} X_{2})
\end{pmatrix*}
 =
\begin{pmatrix*}[l]
Y_{1}	& Y_{2} \\ \addlinespace[1mm]
Y_{2}	& Y_{3}
\end{pmatrix*}.
\end{equation}
It follows that, if $X_{1}$ and $X_{2}$ are not strictly correlated,
\begin{equation*}
\det (\boldsymbol{\Sigma}_{\boldsymbol{X} \boldsymbol{X}})
 = Y_{1} Y_{3} - Y_{2}^{2}.
\end{equation*}
This quantity is strictly positive. Thus,
\begin{equation} \label{eq-13}
\boldsymbol{\Sigma}_{\boldsymbol{X} \boldsymbol{X}}^{-1}
 = \frac{1}{Y_{1} Y_{3} - Y_{2}^{2}}
\begin{pmatrix*}[r]
Y_{3}	& -Y_{2} \\ \addlinespace[1mm]
-Y_{2}	& Y_{1}
\end{pmatrix*}.
\end{equation}

Expanding our definitions, $\boldsymbol{\Sigma}_{\boldsymbol{a} \boldsymbol{a}, \boldsymbol{X}}, \boldsymbol{\Sigma}_{\boldsymbol{b} \boldsymbol{b}, \boldsymbol{X}}, \boldsymbol{\Sigma}_{\boldsymbol{a} \boldsymbol{b}, \boldsymbol{X}}$, and $\boldsymbol{\Sigma}_{\boldsymbol{b} \boldsymbol{a}, \boldsymbol{X}}$ are obtained as follows. Direct evaluation shows that
\begin{equation} \label{eq-14}
\boldsymbol{\Sigma}_{\boldsymbol{a} \boldsymbol{a}}
 = E (a_{j} a_{k})
 = \delta_{j k} \sigma_{a_{j}}^{2}
\end{equation}
and
\begin{equation} \label{eq-15}
\boldsymbol{\Sigma}_{\boldsymbol{b} \boldsymbol{b}}
 = E (b_{j} b_{k})
 = \delta_{j k} \sigma_{b_{j}}^{2},
\end{equation}
where
\begin{equation*}
\delta_{j k}
 =
\left\{ \begin{array}{ll}
1 & \mbox{if $j = k$,} \\
0 & \mbox{if $j \neq k$.}
\end{array}
\right.
\end{equation*}
Further, notice that
\begin{equation} \label{eq-16}
\boldsymbol{\Sigma}_{\boldsymbol{a} \boldsymbol{b}}
 = E (a_{j} b_{k})
 = 0
\end{equation}
and
\begin{equation} \label{eq-17}
\boldsymbol{\Sigma}_{\boldsymbol{b} \boldsymbol{a}}
 = 0.
\end{equation}
Next, since $E (a_{j} X_{1}) = \sigma_{a_{j}}^{2} u_{j}$ and $E (a_{j} X_{2}) = \sigma_{a_{j}}^{2} v_{j}$ for $0 \leq j \leq N$,
\begin{equation} \label{eq-18}
\boldsymbol{\Sigma}_{\boldsymbol{a} \boldsymbol{X}}
 =
(\begin{matrix*}[l]
\sigma_{a_{j}}^{2} u_{j} & \sigma_{a_{j}}^{2} v_{j}
\end{matrix*}),
\end{equation}
whence
\begin{equation} \label{eq-19}
\boldsymbol{\Sigma}_{\boldsymbol{X} \boldsymbol{a}}
 =
\begin{pmatrix*}[l]
\sigma_{a_{k}}^{2} u_{k} \\ \addlinespace[1mm]
\sigma_{a_{k}}^{2} v_{k}
\end{pmatrix*}.
\end{equation}
Analogously, since $E (b_{j} X_{1}) = -\sigma_{b_{j}}^{2} v_{j}$ and $E (b_{j} X_{2}) = \sigma_{b_{j}}^{2} u_{j}$ for $0 \leq j \leq N$,
\begin{equation} \label{eq-20}
\boldsymbol{\Sigma}_{\boldsymbol{b} \boldsymbol{X}}
 =
(\begin{matrix*}[l]
-\sigma_{b_{j}}^{2} v_{j} & \sigma_{b_{j}}^{2} u_{j}
\end{matrix*}),
\end{equation}
whence
\begin{equation} \label{eq-21}
\boldsymbol{\Sigma}_{\boldsymbol{X} \boldsymbol{b}}
 =
\begin{pmatrix*}[r]
-\sigma_{b_{k}}^{2} v_{k} \\ \addlinespace[1mm]
\sigma_{b_{k}}^{2} u_{k}
\end{pmatrix*}
\end{equation}
for $0 \leq k \leq N$.

Then, from \eqref{eq-10}, \eqref{eq-13}, \eqref{eq-14}, \eqref{eq-18} and \eqref{eq-19} for the $j$th row and $k$th column
\begin{equation} \label{eq-22}
\begin{split}
\boldsymbol{\Sigma}_{\boldsymbol{a} \boldsymbol{a}, \boldsymbol{X}}
 & = \delta_{j k} \sigma_{a_{j}}^{2} - \sigma_{a_{j}}^{2} \sigma_{a_{k}}^{2} \left(\frac{Y_{1} v_{j} v_{k} - Y_{2} (u_{j} v_{k} + v_{j} u_{k}) + Y_{3} u_{j} u_{k}}{Y_{1} Y_{3} - Y_{2}^{2}}\right).
\end{split}
\end{equation}
From \eqref{eq-10}, \eqref{eq-13}, \eqref{eq-15}, \eqref{eq-20} and \eqref{eq-21}
\begin{equation} \label{eq-23}
\begin{split}
\boldsymbol{\Sigma}_{\boldsymbol{b} \boldsymbol{b}, \boldsymbol{X}}
 &= \delta_{j k} \sigma_{b_{j}}^{2} - \sigma_{b_{j}}^{2} \sigma_{b_{k}}^{2} \left(\frac{Y_{1} u_{j} u_{k} + Y_{2} (u_{j} v_{k} + v_{j} u_{k}) + Y_{3} v_{j} v_{k}}{Y_{1} Y_{3} - Y_{2}^{2}}\right).
\end{split}
\end{equation}
From \eqref{eq-10}, \eqref{eq-13}, \eqref{eq-16}, \eqref{eq-18} and \eqref{eq-21}
\begin{equation} \label{eq-24}
\boldsymbol{\Sigma}_{\boldsymbol{a} \boldsymbol{b}, \boldsymbol{X}}
 = - \sigma_{a_{j}}^{2} \sigma_{b_{k}}^{2} \left(\frac{Y_{1} v_{j} u_{k} - Y_{2} (u_{j} u_{k} - v_{j} v_{k}) - Y_{3} u_{j} v_{k}}{Y_{1} Y_{3} - Y_{2}^{2}}\right).
\end{equation}
From \eqref{eq-10}, \eqref{eq-13}, \eqref{eq-17}, \eqref{eq-19} and \eqref{eq-20}
\begin{equation} \label{eq-25}
\boldsymbol{\Sigma}_{\boldsymbol{b} \boldsymbol{a}, \boldsymbol{X}}
 = - \sigma_{b_{j}}^{2} \sigma_{a_{k}}^{2} \left(\frac{Y_{1} u_{j} v_{k} - Y_{2} (u_{j} u_{k} - v_{j} v_{k}) - Y_{3} v_{j} u_{k}}{Y_{1} Y_{3} - Y_{2}^{2}}\right).
\end{equation}

The mean function in \eqref{eq-8} is then found by applications of
\begin{equation*}
E (\boldsymbol{a} \mid \boldsymbol{X} = \boldsymbol{K})
 = E (\boldsymbol{a}) + \boldsymbol{\Sigma}_{\boldsymbol{a} \boldsymbol{X}} \boldsymbol{\Sigma}_{\boldsymbol{X} \boldsymbol{X}}^{-1} (\boldsymbol{K} - E (\boldsymbol{X})),
\end{equation*}
which, for the aforesaid reasons, reduces to
\begin{equation} \label{eq-26}
E (\boldsymbol{a} \mid \boldsymbol{X} = \boldsymbol{K})
 = \boldsymbol{\Sigma}_{\boldsymbol{a} \boldsymbol{X}} \boldsymbol{\Sigma}_{\boldsymbol{X} \boldsymbol{X}}^{-1} \boldsymbol{K}.
\end{equation}
From \eqref{eq-7}, \eqref{eq-13}, \eqref{eq-18} and \eqref{eq-26}
\begin{equation} \label{eq-27}
E (a_{j} \mid \boldsymbol{X} = \boldsymbol{K})
 = \sigma_{a_{j}}^{2} \left(\frac{(K_{1} Y_{3} - K_{2} Y_{2}) u_{j} - (K_{1} Y_{2} - K_{2} Y_{1}) v_{j}}{Y_{1} Y_{3} - Y_{2}^{2}}\right).
\end{equation}
The formula for $E (\boldsymbol{b} \mid \boldsymbol{X} = \boldsymbol{K})$ is analogous. Then, from \eqref{eq-7}, \eqref{eq-13}, \eqref{eq-20} and \eqref{eq-26} with the obvious substitution
\begin{equation} \label{eq-28}
E (b_{j} \mid \boldsymbol{X} = \boldsymbol{K})
 = - \sigma_{b_{j}}^{2} \left(\frac{(K_{1} Y_{3} - K_{2} Y_{2}) v_{j} + (K_{1} Y_{2} - K_{2} Y_{1}) u_{j}}{Y_{1} Y_{3} - Y_{2}^{2}}\right).
\end{equation}
Then, from \eqref{eq-9}, \eqref{eq-22}--\eqref{eq-25} and \eqref{eq-27}
\begin{equation} \label{eq-29}
\begin{split}
& E (a_{j} a_{k} \mid \boldsymbol{X} = \boldsymbol{K})
 = E (a_{j} \mid \boldsymbol{X} = \boldsymbol{K}) E (a_{k} \mid \boldsymbol{X} = \boldsymbol{K}) + \operatorname{Cov} (a_{j}, a_{k} \mid \boldsymbol{X} = \boldsymbol{K})
 = \sigma_{a_{j}}^{2} \sigma_{a_{k}}^{2} \\ & \hspace{1em} \times \left(\frac{(K_{1} Y_{3} - K_{2} Y_{2})^{2} u_{j} u_{k} + (K_{1} Y_{2} - K_{2} Y_{1})^{2} v_{j} v_{k} - (K_{1} Y_{3} - K_{2} Y_{2}) (K_{1} Y_{2} - K_{2} Y_{1}) (u_{j} v_{k} + v_{j} u_{k})}{(Y_{1} Y_{3} - Y_{2}^{2})^{2}}\right) \\ & \hspace{2em} + \delta_{j k} \sigma_{a_{j}}^{2} - \sigma_{a_{j}}^{2} \sigma_{a_{k}}^{2} \left(\frac{Y_{3} u_{j} u_{k} + Y_{1} v_{j} v_{k} - Y_{2} (u_{j} v_{k} + v_{j} u_{k})}{Y_{1} Y_{3} - Y_{2}^{2}}\right).
\end{split}
\end{equation}
From \eqref{eq-9}, \eqref{eq-22}--\eqref{eq-25} and \eqref{eq-28}
\begin{equation} \label{eq-30}
\begin{split}
& E (b_{j} b_{k} \mid \boldsymbol{X} = \boldsymbol{K})
 = E (b_{j} \mid \boldsymbol{X} = \boldsymbol{K}) E (b_{k} \mid \boldsymbol{X} = \boldsymbol{K}) + \operatorname{Cov} (b_{j}, b_{k} \mid \boldsymbol{X} = \boldsymbol{K})
 = \sigma_{b_{j}}^{2} \sigma_{b_{k}}^{2} \\ & \hspace{1em} \times \left(\frac{(K_{1} Y_{2} - K_{2} Y_{1})^{2} u_{j} u_{k} + (K_{1} Y_{3} - K_{2} Y_{2})^{2} v_{j} v_{k} + (K_{1} Y_{2} - K_{2} Y_{1}) (K_{1} Y_{3} - K_{2} Y_{2}) (u_{j} v_{k} + v_{j} u_{k})}{(Y_{1} Y_{3} - Y_{2}^{2})^{2}}\right) \\ & \hspace{2em} + \delta_{j k} \sigma_{b_{j}}^{2} - \sigma_{b_{j}}^{2} \sigma_{b_{k}}^{2} \left(\frac{Y_{1} u_{j} u_{k} + Y_{3} v_{j} v_{k} + Y_{2} (u_{j} v_{k} + v_{j} u_{k})}{Y_{1} Y_{3} - Y_{2}^{2}}\right).
\end{split}
\end{equation}
From \eqref{eq-9}, \eqref{eq-22}--\eqref{eq-25}, \eqref{eq-27} and \eqref{eq-28}
\begin{equation} \label{eq-31}
\begin{split}
& E (a_{j} b_{k} \mid \boldsymbol{X} = \boldsymbol{K})
 = E (a_{j} \mid \boldsymbol{X} = \boldsymbol{K}) E (b_{k} \mid \boldsymbol{X} = \boldsymbol{K}) + \operatorname{Cov} (a_{j}, b_{k} \mid \boldsymbol{X} = \boldsymbol{K})
 = \sigma_{a_{j}}^{2} \sigma_{b_{k}}^{2} \\ & \hspace{1em} \times \left(\frac{(K_{1} Y_{2} - K_{2} Y_{1})^{2} v_{j} u_{k} - (K_{1} Y_{3} - K_{2} Y_{2})^{2} u_{j} v_{k} - (K_{1} Y_{3} - K_{2} Y_{2}) (K_{1} Y_{2} - K_{2} Y_{1}) (u_{j} u_{k} - v_{j} v_{k})}{(Y_{1} Y_{3} - Y_{2}^{2})^{2}}\right) \\ & \hspace{2em} - \sigma_{a_{j}}^{2} \sigma_{b_{k}}^{2} \left(\frac{Y_{1} v_{j} u_{k} - Y_{2} (u_{j} u_{k} - v_{j} v_{k}) - Y_{3} u_{j} v_{k}}{Y_{1} Y_{3} - Y_{2}^{2}}\right)
\end{split}
\end{equation}
and, likewise,
\begin{equation} \label{eq-32}
\begin{split}
& E (b_{j} a_{k} \mid \boldsymbol{X} = \boldsymbol{K})
 = E (b_{j} \mid \boldsymbol{X} = \boldsymbol{K}) E (a_{k} \mid \boldsymbol{X} = \boldsymbol{K}) + \operatorname{Cov} (b_{j}, a_{k} \mid \boldsymbol{X} = \boldsymbol{K})
 = \sigma_{b_{j}}^{2} \sigma_{a_{k}}^{2} \\ & \hspace{1em} \times \left(\frac{(K_{1} Y_{2} - K_{2} Y_{1})^{2} u_{j} v_{k} - (K_{1} Y_{3} - K_{2} Y_{2})^{2} v_{j} u_{k} - (K_{1} Y_{3} - K_{2} Y_{2}) (K_{1} Y_{2} - K_{2} Y_{1}) (u_{j} u_{k} - v_{j} v_{k})}{(Y_{1} Y_{3} - Y_{2}^{2})^{2}}\right) \\ & \hspace{2em} - \sigma_{b_{j}}^{2} \sigma_{a_{k}}^{2} \left(\frac{Y_{1} u_{j} v_{k} - Y_{2} (u_{j} u_{k} - v_{j} v_{k}) - Y_{3} v_{j} u_{k}}{Y_{1} Y_{3} - Y_{2}^{2}}\right).
\end{split}
\end{equation}
Then, from \eqref{eq-29} and \eqref{eq-30}
\begin{equation} \label{eq-33}
\begin{split}
& E (a_{j} a_{k} + b_{j} b_{k} \mid \boldsymbol{X} = \boldsymbol{K})
 = \frac{1}{(Y_{1} Y_{3} - Y_{2}^{2})^{2}} ((K_{1} Y_{3} - K_{2} Y_{2})^{2} (\sigma_{a_{j}}^{2} \sigma_{a_{k}}^{2} u_{j} u_{k} + \sigma_{b_{j}}^{2} \sigma_{b_{k}}^{2} v_{j} v_{k}) \\ & \hspace{5em} + (K_{1} Y_{2} - K_{2} Y_{1})^{2} (\sigma_{a_{j}}^{2} \sigma_{a_{k}}^{2} v_{j} v_{k} + \sigma_{b_{j}}^{2} \sigma_{b_{k}}^{2} u_{j} u_{k}) \\ & \hspace{5em} - (K_{1} Y_{3} - K_{2} Y_{2}) (K_{1} Y_{2} - K_{2} Y_{1}) (u_{j} v_{k} + v_{j} u_{k}) (\sigma_{a_{j}}^{2} \sigma_{a_{k}}^{2} - \sigma_{b_{j}}^{2} \sigma_{b_{k}}^{2})) \\ & \hspace{5em} - \frac{1}{Y_{1} Y_{3} - Y_{2}^{2}} (Y_{1} (\sigma_{a_{j}}^{2} \sigma_{a_{k}}^{2} v_{j} v_{k} + \sigma_{b_{j}}^{2} \sigma_{b_{k}}^{2} u_{j} u_{k}) + Y_{2} (u_{j} v_{k} + v_{j} u_{k}) (\sigma_{b_{j}}^{2} \sigma_{b_{k}}^{2} - \sigma_{a_{j}}^{2} \sigma_{a_{k}}^{2}) \\ & \hspace{5em} + Y_{3} (\sigma_{a_{j}}^{2} \sigma_{a_{k}}^{2} u_{j} u_{k} + \sigma_{b_{j}}^{2} \sigma_{b_{k}}^{2} v_{j} v_{k})) + \delta_{j k} (\sigma_{a_{j}}^{2} + \sigma_{b_{j}}^{2}).
	\end{split}
\end{equation}
From \eqref{eq-31} and \eqref{eq-32}
\begin{equation} \label{eq-34}
\begin{split}
& E (a_{j} b_{k} - b_{j} a_{k} \mid \boldsymbol{X}= \boldsymbol{K})
 = \frac{1}{(Y_{1} Y_{3} - Y_{2}^{2})^{2}} ((K_{1} Y_{2} - K_{2} Y_{1})^{2} (\sigma_{a_{j}}^{2} \sigma_{b_{k}}^{2} v_{j} u_{k} - \sigma_{b_{j}}^{2} \sigma_{a_{k}}^{2} u_{j} v_{k}) \\ & \hspace{5em} - (K_{1} Y_{3} - K_{2} Y_{2})^{2} (\sigma_{a_{j}}^{2} \sigma_{b_{k}}^{2} u_{j} v_{k} - \sigma_{b_{j}}^{2} \sigma_{a_{k}}^{2} v_{j} u_{k}) \\ & \hspace{5em} - (K_{1} Y_{3} - K_{2} Y_{2}) (K_{1} Y_{2} - K_{2} Y_{1}) (u_{j} u_{k} - v_{j} v_{k}) (\sigma_{a_{j}}^{2} \sigma_{b_{k}}^{2} - \sigma_{b_{j}}^{2} \sigma_{a_{k}}^{2})) \\ & \hspace{5em} - \frac{1}{Y_{1} Y_{3} - Y_{2}^{2}} (Y_{1} (\sigma_{a_{j}}^{2} \sigma_{b_{k}}^{2} v_{j} u_{k} - \sigma_{b_{j}}^{2} \sigma_{a_{k}}^{2} u_{j} v_{k}) + Y_{2} (v_{j} v_{k} - u_{j} u_{k}) (\sigma_{a_{j}}^{2} \sigma_{b_{k}}^{2} - \sigma_{b_{j}}^{2} \sigma_{a_{k}}^{2}) \\ & \hspace{5em} - Y_{3} (\sigma_{a_{j}}^{2} \sigma_{b_{k}}^{2} u_{j} v_{k} - \sigma_{b_{j}}^{2} \sigma_{a_{k}}^{2} v_{j} u_{k})).
\end{split}
\end{equation}

Altogether, in view of \eqref{eq-5}, \eqref{eq-33} and \eqref{eq-34}, after all the necessary simplifications,
\begin{equation} \label{eq-35}
\begin{split}
& E (\det (\nabla \boldsymbol{X}) \mid \boldsymbol{X} = \boldsymbol{K})
 = \sum_{j = 0}^{N} (\sigma_{a_{j}}^{2} + \sigma_{b_{j}}^{2}) \left(\left(\frac{\partial u_{j}}{\partial x}\right)^{2} + \left(\frac{\partial v_{j}}{\partial x}\right)^{2}\right) \\ & \hspace{4em} - \left(\frac{Y_{1} + Y_{2}}{Y_{1} Y_{3} - Y_{2}^{2}} - \frac{(K_{1} Y_{2} - K_{2} Y_{1})^{2} + (K_{1} Y_{3} - K_{2} Y_{2}) (K_{1} Y_{2} - K_{2} Y_{1})}{(Y_{1} Y_{3} - Y_{2}^{2})^{2}}\right) \\ & \hspace{12em} \times \left|\sum_{j = 0}^{N} (\sigma_{b_{j}}^{2} u_{j} - i \sigma_{a_{j}}^{2} v_{j}) \left(\frac{\partial u_{j}}{\partial x} + i \frac{\partial v_{j}}{\partial x}\right)\right|^{2} \\ & \hspace{4em} - \left(\frac{Y_{2} + Y_{3}}{Y_{1} Y_{3} - Y_{2}^{2}} - \frac{(K_{1} Y_{3} - K_{2} Y_{2})^{2} + (K_{1} Y_{3} - K_{2} Y_{2}) (K_{1} Y_{2} - K_{2} Y_{1})}{(Y_{1} Y_{3} - Y_{2}^{2})^{2}}\right) \\ & \hspace{12em} \times \left|\sum_{j = 0}^{N} (\sigma_{a_{j}}^{2} u_{j} - i \sigma_{b_{j}}^{2} v_{j}) \left(\frac{\partial u_{j}}{\partial x} + i \frac{\partial v_{j}}{\partial x}\right)\right|^{2} \\ & \hspace{4em} + \left(\frac{Y_{2}}{Y_{1} Y_{3} - Y_{2}^{2}} - \frac{(K_{1} Y_{3} - K_{2} Y_{2}) (K_{1} Y_{2} - K_{2} Y_{1})}{(Y_{1} Y_{3} - Y_{2}^{2})^{2}}\right) \\ & \hspace{12em} \times \left|\sum_{j = 0}^{N} ((\sigma_{a_{j}}^{2} u_{j} - i \sigma_{b_{j}}^{2} v_{j}) + i (\sigma_{b_{j}}^{2} u_{j} - i \sigma_{a_{j}}^{2} v_{j})) \left(\frac{\partial u_{j}}{\partial x} + i \frac{\partial v_{j}}{\partial x}\right)\right|^{2}.
\end{split}
\end{equation}
Since $X_{1}$ and $X_{2}$ are random variables distributed according to the normal law, their joint density is
\begin{equation} \label{eq-36}
\begin{split}
p_{X_{1}, X_{2}} (\boldsymbol{K}^{\prime})
 & = \frac{1}{2 \pi \sqrt{\det (\boldsymbol{\Sigma}_{\boldsymbol{X} \boldsymbol{X}})}} \exp \left(- \frac{1}{2} (\boldsymbol{K} - E (\boldsymbol{X}))^{\prime} \boldsymbol{\Sigma}_{\boldsymbol{X} \boldsymbol{X}}^{-1} (\boldsymbol{K} - E (\boldsymbol{X}))\right) \\
 & = \frac{1}{2 \pi \sqrt{Y_{1} Y_{3} - Y_{2}^{2}}} \exp \left(-\frac{K_{1}^{2} Y_{3} - 2 K_{1} K_{2} Y_{2} + K_{2}^{2} Y_{1}}{2 (Y_{1} Y_{3} - Y_{2}^{2})}\right).
\end{split}
\end{equation}
Hence, in accordance with \eqref{eq-8}, \eqref{eq-35} and \eqref{eq-36}, the required result is proved.

\section{Reformulation of Theorem \ref{th-2}} \label{sec-3}

If $\{A_{j}\}_{j = 0}^{\infty}$ and $\{B_{j}\}_{j = 0}^{\infty}$ are sequences of i.i.d. real normal random variables for which the respective increments $A_{j} - A_{j - 1}$ and $B_{j} - B_{j - 1}$ are independent for $j \geq 0$ and $A_{-1} = B_{-1} = 0$ by convention, then the increments 
\begin{equation*}
\Delta_{j}
 = (A_{j} - A_{j - 1}) + i (B_{j} - B_{j - 1})
\end{equation*}
are independent complex normal random variables with mean zero and finite $\operatorname{Var} (\Delta_{j})$ such that $A_{j} + i B_{j} = \Delta_{0} + \Delta_{1} + \cdots + \Delta_{j}$ for $j \geq 0$. Then $\{A_{j} + i B_{j}\}_{j = 0}^{\infty}$ can be interpreted as a sequence of successive observations of a Brownian motion. More precisely, $A_{j} + i B_{j} = W (t_{j})$ for $j \geq 0$, where $t_{0} < t_{1} < \ldots$ and $\{W (t)\}_{t = 0}^{\infty}$ is the standard Brownian motion. It is plain that $\operatorname{Var} (\Delta_{j})$ is the distance between the successive times $t_{j - 1}$ and $t_{j}$ for $j \geq 0$. Thus, the sum in \eqref{eq-2} assumes the form
\begin{equation} \label{eq-37}
S_{N} (z)
 = \sum_{j = 0}^{N} (A_{j} + i B_{j}) f_{j} (z)
 = \sum_{k = 0}^{N} F_{k} (z) \Delta_{k},
\end{equation}
where
\begin{equation} \label{eq-38}
F_{k} (z)
 = \sum_{j = k}^{N} u_{j} (x, y) + i \sum_{j = k}^{N} v_{j} (x, y)
\end{equation}
for $0 \leq k \leq N$. Observe that $\{F_{k} (z)\}_{k = 0}^{N}$ is a sequence of holomorphic functions that are real-valued on $\mathds{R}$. Hence, $\overline{F_{k} (z)} = F_{k} (\overline{z})$ for $0 \leq k \leq N$ and all $z \in \mathds{C}$. The covariance matrix of $\Delta_{k}$ is given by
\begin{equation*}
\Gamma_{k}
 = 
\begin{pmatrix*}[c]
\sigma_{a_{k}}^{2} & 0 \\ \addlinespace[1mm]
0 & \sigma_{b_{k}}^{2}
\end{pmatrix*}
\end{equation*}
for $0 \leq k \leq N$. Then, from Theorem \ref{th-2} the following result is attained.

\begin{theorem}
Provided all the conditions imposed on $S_{N} (z)$ in \eqref{eq-37} and \eqref{eq-38} and $T$ are satisfied, then for all $N > 1$ the formula for $h_{N, \boldsymbol{K}} (z)$ in Theorem \ref{th-2} now holds for
\begin{equation*}
	D_{0} (z)
	= \sqrt{Y_{1} (z) Y_{3} (z) - Y_{2}^{2} (z)},
\end{equation*}
where
\begin{equation*}
\begin{array}{c@{\qquad}c}
\displaystyle Y_{1} (z)
 = \sum_{k = 0}^{N} \left(\sigma_{a_{k}}^{2} \left(\sum_{j = k}^{N} u_{j}\right)^{2} + \sigma_{b_{k}}^{2} \left(\sum_{j = k}^{N} v_{j}\right)^{2}\right), \\ [5ex]
\displaystyle Y_{2} (z)
 = \sum_{k = 0}^{N} \left((\sigma_{a_{k}}^{2} - \sigma_{b_{k}}^{2}) \left(\sum_{j = k}^{N} u_{j}\right) \left(\sum_{j = k}^{N} v_{j}\right)\right), \\ [5ex]
\displaystyle Y_{3} (z)
 = \sum_{k = 0}^{N} \left(\sigma_{b_{k}}^{2} \left(\sum_{j = k}^{N} u_{j}\right)^{2} + \sigma_{a_{k}}^{2} \left(\sum_{j = k}^{N} v_{j}\right)^{2}\right),
\end{array}
\end{equation*}
and
\begin{equation*}
\begin{array}{c@{\qquad}c}
\displaystyle D_{1} (z)
 = \sum_{k = 0}^{N} \left(\sigma_{a_{k}}^{2} \sum_{j = k}^{N} u_{j} - i \sigma_{b_{k}}^{2} \sum_{j = k}^{N} v_{j}\right) \left(\sum_{j = k}^{N} \frac{\partial u_{j}}{\partial x} + i \sum_{j = k}^{N} \frac{\partial v_{j}}{\partial x}\right), \\ [5ex]
\displaystyle D_{2} (z)
 = \sum_{k = 0}^{N} \left(\sigma_{b_{k}}^{2} \sum_{j = k}^{N} u_{j} - i \sigma_{a_{k}}^{2} \sum_{j = k}^{N} v_{j}\right) \left(\sum_{j = k}^{N} \frac{\partial u_{j}}{\partial x} + i \sum_{j = k}^{N} \frac{\partial v_{j}}{\partial x}\right), \\ [5ex]
\displaystyle D_{3} (z)
 = \sum_{k = 0}^{N} (\sigma_{a_{k}}^{2} + \sigma_{b_{k}}^{2}) \left(\left(\sum_{j = k}^{N} \frac{\partial u_{j}}{\partial x}\right)^{2} + \left(\sum_{j = k}^{N} \frac{\partial v_{j}}{\partial x}\right)^{2}\right).
\end{array}
\end{equation*}
\end{theorem}

\section{Proof of Theorem \ref{th-4}} \label{sec-4}

The proof of Theorem \ref{th-4} mirrors that of Theorem \ref{th-2}. Below, only the differences are highlighted. Direct computation leads to the equalities
\begin{equation*}
\begin{array}{c@{\qquad}c}
\displaystyle E (X_{1})
 = \sum_{j = 0}^{N} (\mu_{a_{j}} u_{j} - \mu_{b_{j}} v_{j}), \qquad
\displaystyle E (X_{2})
 = \sum_{j = 0}^{N} (\mu_{a_{j}} v_{j} + \mu_{b_{j}} u_{j}), \\ [3ex]
\displaystyle E (X_{1}^{2})
= \sum_{j = 0}^{N} (\sigma_{a_{j}}^{2} u_{j}^{2} + \sigma_{b_{j}}^{2} v_{j}^{2}), \qquad
\displaystyle E (X_{2}^{2})
= \sum_{j= 0}^{N} (\sigma_{a_{j}}^{2} v_{j}^{2} + \sigma_{b_{j}}^{2} u_{j}^{2}).
\end{array}
\end{equation*}
Further,
\begin{equation*}
\displaystyle E (X_{1} X_{2})
 = \sum_{j = 0}^{N} (\sigma_{a_{j}}^{2} - \sigma_{b_{j}}^{2}) u_{j} v_{j}.
\end{equation*}
It follows that
\begin{equation*}
\boldsymbol{\Sigma}_{\boldsymbol{X} \boldsymbol{X}}
 =
\begin{pmatrix*}[l]
Y_{1}^{\ast}	& Y_{2}^{\ast} \\ \addlinespace[1mm]
Y_{2}^{\ast}	& Y_{3}^{\ast}
\end{pmatrix*},
\end{equation*}
and hence
\begin{equation*}
\boldsymbol{\Sigma}_{\boldsymbol{X} \boldsymbol{X}}^{-1}
 = \frac{1}{Y_{1}^{\ast} Y_{3}^{\ast} - (Y_{2}^{\ast})^{2}}
\begin{pmatrix*}[r]
Y_{3}^{\ast}	& -Y_{2}^{\ast} \\ \addlinespace[1mm]
-Y_{2}^{\ast}	& Y_{1}^{\ast}
\end{pmatrix*}.
\end{equation*}

On expanding our definitions, for the $j$th row and $k$th column
\begin{equation*}
\begin{split}
\boldsymbol{\Sigma}_{\boldsymbol{a} \boldsymbol{a}}
 & = \delta_{j k} \sigma_{a_{j}}^{2} - \mu_{a_{j}} \mu_{a_{k}}, \\
\boldsymbol{\Sigma}_{\boldsymbol{b} \boldsymbol{a}}
& = - \mu_{b_{j}} \mu_{a_{k}},
\end{split}
\qquad
\begin{split}
\boldsymbol{\Sigma}_{\boldsymbol{a} \boldsymbol{b}}
& = - \mu_{a_{j}} \mu_{b_{k}}, \\
\boldsymbol{\Sigma}_{\boldsymbol{b} \boldsymbol{b}}
& = \delta_{j k} \sigma_{b_{j}}^{2} - \mu_{b_{j}} \mu_{b_{k}}.
\end{split}
\end{equation*}
Since
\begin{equation*}
\begin{split}
E ((a_{j} - E (a_{j})) (X_{1} - E (X_{1})))
& = A_{j, 1}, \\
E ((b_{j} - E (b_{j})) (X_{1} - E (X_{1})))
& = - B_{j, 1},
\end{split}
\qquad
\begin{split}
E ((a_{j} - E (a_{j})) (X_{2} - E (X_{2})))
& = A_{j, 2}, \\
E ((b_{j} - E (b_{j})) (X_{2} - E (X_{2})))
& = B_{j, 2},
\end{split}
\end{equation*}
we have
\begin{equation*}
\boldsymbol{\Sigma}_{\boldsymbol{a} \boldsymbol{X}}
 =
(\begin{matrix*}[l]
A_{j, 1} & A_{j, 2}
\end{matrix*}), \qquad
\boldsymbol{\Sigma}_{\boldsymbol{a} \boldsymbol{X}}
 =
(\begin{matrix*}[l]
- B_{j, 1} & B_{j, 2}
\end{matrix*}).
\end{equation*}Using \eqref{eq-10}, simple algebra leads to
\begin{align*}
\boldsymbol{\Sigma}_{\boldsymbol{a} \boldsymbol{a}, \boldsymbol{X}}
& = \delta_{j k} \sigma_{a_{j}}^{2} - \mu_{a_{j}} \mu_{a_{k}} - \frac{A_{j, 1} A_{k, 1} Y_{3}^{\ast} + A_{j, 2} A_{k, 2} Y_{1}^{\ast} - (A_{j, 1} A_{k, 2} + A_{j, 2} A_{k, 1}) Y_{2}^{\ast}}{Y_{1}^{\ast} Y_{3}^{\ast} - (Y_{2}^{\ast})^{2}}, \\
\boldsymbol{\Sigma}_{\boldsymbol{b} \boldsymbol{b}, \boldsymbol{X}}
& = \delta_{j k} \sigma_{b_{j}}^{2} - \mu_{b_{j}} \mu_{b_{k}} - \frac{B_{j, 1} B_{k, 1} Y_{3}^{\ast} + B_{j, 1} B_{k, 2} Y_{1}^{\ast} + (B_{j, 1} B_{k, 2} + B_{j, 2} B_{k, 1}) Y_{2}^{\ast}}{Y_{1}^{\ast} Y_{3}^{\ast} - (Y_{2}^{\ast})^{2}}, \\
\boldsymbol{\Sigma}_{\boldsymbol{a} \boldsymbol{b}, \boldsymbol{X}}
& = - \mu_{a_{j}} \mu_{b_{k}} + \frac{A_{j, 1} B_{k, 1} Y_{3}^{\ast} - A_{j, 2} B_{k, 2} Y_{1}^{\ast} + (A_{j, 1} B_{k, 2} - A_{j, 2} B_{k, 1}) Y_{2}^{\ast}}{Y_{1}^{\ast} Y_{3}^{\ast} - (Y_{2}^{\ast})^{2}}, \\
\boldsymbol{\Sigma}_{\boldsymbol{b} \boldsymbol{a}, \boldsymbol{X}}
& = - \mu_{b_{j}} \mu_{a_{k}} + \frac{B_{j, 1} A_{k, 1} Y_{3}^{\ast} - B_{j, 2} A_{k, 2} Y_{1}^{\ast} + (B_{j, 2} A_{k, 1} - B_{j, 1} A_{k, 1}) Y_{2}^{\ast}}{Y_{1}^{\ast} Y_{3}^{\ast} - (Y_{2}^{\ast})^{2}}.
\end{align*}
Using these in \eqref{eq-26}, we obtain
\begin{equation*}
\begin{split}
& E (\boldsymbol{a} \mid \boldsymbol{X} = \boldsymbol{K})
 = \mu_{a_{j}} \\ & \hspace{1em} - \frac{((K_{1} - E (X_{1})) Y_{3}^{\ast} - (K_{2} - E (X_{2})) Y_{2}^{\ast}) A_{j, 1} + ((K_{2} - E (X_{2})) Y_{1}^{\ast} - (K_{1} - E (X_{1})) Y_{2}^{\ast}) A_{j, 2}}{Y_{1}^{\ast} Y_{3}^{\ast} - (Y_{2}^{\ast})^{2}}
\end{split}
\end{equation*}
and, likewise,
\begin{equation*}
\begin{split}
& E (\boldsymbol{b} \mid \boldsymbol{X} = \boldsymbol{K})
 = \mu_{b_{j}} \\ & \hspace{1em} + \frac{((K_{1} - E (X_{1})) Y_{3}^{\ast} - (K_{2} - E (X_{2})) Y_{2}^{\ast}) B_{j, 1} - ((K_{2} - E (X_{2})) Y_{1}^{\ast} - (K_{1} - E (X_{1})) Y_{2}^{\ast}) B_{j, 2}}{Y_{1}^{\ast} Y_{3}^{\ast} - (Y_{2}^{\ast})^{2}}.
\end{split}
\end{equation*}
Thus, for the $j$th row and $k$th column
\begingroup
\allowdisplaybreaks
\begin{align*}
& E (a_{j} a_{k} \mid \boldsymbol{X} = \boldsymbol{K})
 = \delta_{j k} \sigma_{a_{j}}^{2} - (A_{j, 1} A_{k, 1} Y_{3}^{\ast} + A_{j, 2} A_{k, 2} Y_{1}^{\ast} - (A_{j, 1} A_{k, 2} + A_{j, 2} A_{k, 1}) Y_{2}^{\ast}) / (Y_{1}^{\ast} Y_{3}^{\ast} - (Y_{2}^{\ast})^{2}) \\ & \hspace{1em} - (((K_{1} - E (X_{1})) Y_{3}^{\ast} - (K_{2} - E (X_{2})) Y_{2}^{\ast}) (\mu_{a_{j}} A_{k, 1} + \mu_{a_{k}} A_{j, 1}) \\ & \hspace{1em} + ((K_{2} - E (X_{2})) Y_{1}^{\ast} - (K_{1} - E (X_{1})) Y_{2}^{\ast}) (\mu_{a_{j}} A_{k, 2} + \mu_{a_{k}} A_{j, 2}) / (Y_{1}^{\ast} Y_{3}^{\ast} - (Y_{2}^{\ast})^{2}) \\ & \hspace{1em} + ((K_{1} - E (X_{1})) Y_{3}^{\ast} - (K_{2} - E (K_{2})) Y_{2}^{\ast}) ((K_{2} - E (X_{2})) Y_{1}^{\ast} - (K_{1} - E (X_{1})) (Y_{2}^{\ast})^{\ast}) \\ &\hspace{1em} \times (A_{j, 1} A_{k, 2} + A_{j, 2} A_{k, 1}) / (Y_{1}^{\ast} Y_{3}^{\ast} - (Y_{2}^{\ast})^{2})^{2} \\ & \hspace{1em} + (((K_{1} + E (X_{1})) Y_{3}^{\ast} - (K_{2} - E (X_{2})) Y_{2}^{\ast})^{2} A_{j, 1} A_{k, 1} \\ & \hspace{1em} + ((K_{2} - E (X_{2})) Y_{1}^{\ast} - (K_{1} - E (X_{1})) (Y_{2}^{\ast})^{2} A_{j, 2} A_{k, 2}) / (Y_{1}^{\ast} Y_{3}^{\ast} - (Y_{2}^{\ast})^{2})^{2}, \\ \\
& E (b_{j} b_{k} \mid \boldsymbol{X} = \boldsymbol{K})
 = \delta_{j k} \sigma_{b_{j}}^{2} - (B_{j, 1} B_{k, 1} Y_{3}^{\ast} + B_{j, 2} B_{k, 2} Y_{1}^{\ast} + (B_{j, 1} B_{k, 2} + B_{j, 2} B_{k, 1}) Y_{2}^{\ast}) / (Y_{1}^{\ast} Y_{3}^{\ast} - (Y_{2}^{\ast})^{2}) \\ &\hspace{1em} + ((K_{1} - E (X_{1})) Y_{3}^{\ast} - (K_{2} - E (X_{2})) Y_{2}^{\ast}) (\mu_{b_{j}} B_{k, 1} + \mu_{b_{k}} B_{j, 1}) \\ & \hspace{1em} - ((K_{2} - E (X_{2})) Y_{1}^{\ast} - (K_{1} - E (X_{1})) Y_{2}^{\ast}) (\mu_{b_{j}} B_{k, 2} + \mu_{b_{k}} B_{j, 2}) / (Y_{1}^{\ast} Y_{3}^{\ast} - (Y_{2}^{\ast})^{2}) \\ & \hspace{1em} - ((K_{1} - E (X_{1})) Y_{3}^{\ast} - (K_{2} - E (X_{2})) Y_{2}^{\ast}) ((K_{2} - E (X_{2})) Y_{1}^{\ast} - (K_{1} - E (X_{1})) Y_{2}^{\ast}) \\ & \hspace{1em} \times (B_{j, 1} B_{k, 2} + B_{j, 2} B_{k, 1}) / (Y_{1}^{\ast} Y_{3}^{\ast} - (Y_{2}^{\ast})^{2})^{2} \\ & \hspace{1em} + (((K_{1} - E (X_{1})) Y_{3}^{\ast} - (K_{2} - E (X_{2}) Y_{2}^{\ast})^{2} B_{j, 1} B_{k, 1} \\ & \hspace{1em} + ((K_{2} - E (X_{2})) Y_{1}^{\ast} - (K_{1} - E (X_{1})) Y_{2}^{\ast})^{2}) B_{j, 2} B_{k, 2}) / (Y_{1}^{\ast} Y_{3}^{\ast} - (Y_{2}^{\ast})^{2})^{2}, \\ \\
& E (a_{j} b_{k} \mid \boldsymbol{X} = \boldsymbol{K})
 = (A_{j, 1} B_{k, 1} Y_{3}^{\ast} - A_{j, 2} B_{k, 2} Y_{1}^{\ast} + (A_{j, 1} B_{k, 2} - A_{j, 2} B_{k, 1}) Y_{2}^{\ast}) / (Y_{1}^{\ast} Y_{3}^{\ast} - (Y_{2}^{\ast})^{2}) \\ & \hspace{1em} - (((K_{1} - E (X_{1})) Y_{3}^{\ast} - (K_{2} - E (X_{2})) Y_{2}^{\ast}) (\mu_{b_{k}} A_{j, 1} - \mu_{a_{j}} B_{k, 1}) \\ & \hspace{1em} + ((K_{2} - E (X_{2})) Y_{1}^{\ast} - (K_{1} - E (X_{1})) Y_{2}^{\ast}) (\mu_{a_{j}} B_{k, 2} + \mu_{b_{k}} A_{j, 2})) / (Y_{1}^{\ast} Y_{3}^{\ast} - (Y_{2}^{\ast})^{2}) \\ & \hspace{1em} + ((K_{1} - E (X_{1})) Y_{3}^{\ast} - (K_{2} - E (X_{2})) Y_{2}^{\ast}) ((K_{2} - E (X_{2})) Y_{1}^{\ast} - (K_{1} - E (X_{1}) ) Y_{2}^{\ast}) \\ & \hspace{1em} \times (A_{j, 1} B_{k, 2} - A_{j, 2} B_{k, 1}) / (Y_{1}^{\ast} Y_{3}^{\ast} - (Y_{2}^{\ast})^{2})^{2} \\ & \hspace{1em} + (((K_{1} - E (X_{1})) Y_{3}^{\ast} - (K_{2} - E (X_{2})) Y_{2}^{\ast})^{2} A_{j, 1} B_{k, 1} \\ & \hspace{1em} + ((K_{2} - E (X_{2})) Y_{1}^{\ast} - (K_{1} - E (X_{1})) Y_{2}^{\ast})^{2} A_{j, 2} B_{k, 2}) / (Y_{1}^{\ast} Y_{3}^{\ast} - (Y_{2}^{\ast})^{2})^{2}, \\ \\
& E (b_{j} a_{k} \mid \boldsymbol{X} = \boldsymbol{K})
 = (B_{j, 1} A_{k, 1} Y_{3}^{\ast} - B_{j, 2} A_{k, 2} Y_{1}^{\ast} + (B_{j, 2} A_{k, 1} - B_{j, 1} A_{k, 2}) Y_{2}^{\ast}) / (Y_{1}^{\ast} Y_{3}^{\ast} - (Y_{2}^{\ast})^{2}) \\ & \hspace{1em} - (((K_{1} - E (X_{1})) Y_{3}^{\ast} - (K_{2} - E (X_{2})) Y_{2}^{\ast}) (\mu_{b_{j}} A_{k, 1} - \mu_{a_{k}} B_{j, 1}) \\ & \hspace{1em} + ((K_{2} - E (X_{2})) Y_{1}^{\ast} - (K_{1} - E (X_{1})) Y_{2}^{\ast}) (\mu_{b_{j}} A_{k, 2} + \mu_{a_{k}} B_{j, 2})) / (Y_{1}^{\ast} Y_{3}^{\ast} - (Y_{2}^{\ast})^{2}) \\ & \hspace{1em} + ((K_{1} - E (X_{1})) Y_{3}^{\ast} - (K_{2} - E (X_{2})) Y_{2}^{\ast}) ((K_{2} - E (X_{2})) Y_{1}^{\ast} - (K_{1} - E (X_{1}) ) Y_{2}^{\ast}) \\ & \hspace{1em} \times (B_{j, 2} A_{k, 1} - B_{j, 1} A_{k, 2}) / (Y_{1}^{\ast} Y_{3}^{\ast} - (Y_{2}^{\ast})^{2})^{2} \\ & \hspace{1em} - (((K_{1} - E (X_{1})) Y_{3}^{\ast} - (K_{2} - E (X_{2})) Y_{2}^{\ast})^{2} B_{j, 1} A_{k, 1} \\ & \hspace{1em} - ((K_{2} - E (X_{2})) Y_{1}^{\ast} - (K_{1} - E (X_{1})) Y_{2}^{\ast})^{2} B_{j, 2} A_{k, 2}) / (Y_{1}^{\ast} Y_{3}^{\ast} - (Y_{2}^{\ast})^{2})^{2}.
\end{align*}
\endgroup
Having obtained the four expectations above, the required expectations for computing the value of $E (\det (\nabla \boldsymbol{X}) \mid \boldsymbol{X} = \boldsymbol{K})$ can now be derived. Thus, for $j \neq k$
\begingroup
\allowdisplaybreaks
\begin{align*}
& E (a_{j} a_{k} + b_{j} b_{k} \mid \boldsymbol{X} = \boldsymbol{K})
 = (((K_{1} - E (X_{1})) Y_{3}^{\ast} - (K_{2} - E (X_{2})) Y_{2}^{\ast})^{2} (A_{j, 1} A_{k, 1} + B_{j, 1} B_{k, 1}) \\ & \hspace{1em} + ((K_{1} - E (X_{1})) Y_{1}^{\ast} - (K_{2} - E (X_{2})) Y_{2}^{\ast})^{2} (A_{j, 2} A_{k, 2} + B_{j, 2} B_{k, 2})) / (Y_{1}^{\ast} Y_{3}^{\ast} - (Y_{2}^{\ast})^{2}) \\ & \hspace{1em} - (((K_{1} - E (X_{1})) Y_{3}^{\ast} - (K_{2} - E (X_{2})) Y_{2}^{\ast}) (\mu_{a_{j}} A_{k, 1} + \mu_{a_{k}} A_{j, 1} - \mu_{b_{j}} B_{k, 1} - \mu_{b_{k}} B_{j, 1}) \\ & \hspace{1em} + ((K_{2} - E (X_{2})) Y_{1}^{\ast} - (K_{1} - E (X_{1})) Y_{2}^{\ast}) \\ & \hspace{1em} \times (\mu_{a_{j}} A_{k, 2} + \mu_{a_{k}} A_{j, 2} + \mu_{b_{j}} B_{k, 2} + \mu_{b_{k}} B_{j, 2})) / (Y_{1}^{\ast} Y_{3}^{\ast} - (Y_{2}^{\ast})^{2}) \\ & \hspace{1em} + (((K_{1} - E (X_{1})) Y_{3}^{\ast} - (K_{2} - E (X_{2})) Y_{2}^{\ast}) ((K_{2} - E (X_{2})) Y_{1}^{\ast} - (K_{1} - E (X_{1})) Y_{2}^{\ast}) \\ & \hspace{1em} \times (A_{j, 1} A_{k, 2} + A_{j, 2} A_{k, 1} - B_{j, 1} B_{k, 2} - B_{j, 2} B_{k, 1})) / (Y_{1}^{\ast} Y_{3}^{\ast} - (Y_{2}^{\ast})^{2})^{2} \\ & \hspace{1em} + (((K_{1} - E (X_{1})) Y_{3}^{\ast} - (K_{2} - E (X_{2})) Y_{2}^{\ast})^{2} (A_{j, 1} A_{k, 1} + B_{j, 1} B_{k, 1}) \\ & \hspace{1em} + ((K_{2} - E (X_{2})) Y_{1}^{\ast} - (K_{1} - E (X_{1})) Y_{2}^{\ast})^{2} (A_{j, 2} A_{k, 2} + B_{j, 2} B_{k, 2})) / (Y_{1}^{\ast} Y_{3}^{\ast} - (Y_{2}^{\ast})^{2})^{2}, \\ \\
& E (a_{j} b_{k} - b_{j} a_{k} \mid \boldsymbol{X} = \boldsymbol{K})
 = - ((B_{j, 1} A_{k, 1} - A_{j, 1} B_{k, 1}) Y_{3}^{\ast} + (A_{j, 2} B_{k, 2} - B_{j, 2} A_{k, 2}) Y_{1}^{\ast} \\ & \hspace{1em} - (A_{j, 1} B_{k, 2} + B_{j, 1} A_{k, 2} - A_{j, 2} B_{k, 1} - B_{j, 2} A_{k, 1}) Y_{2}^{\ast}) / (Y_{1}^{\ast} Y_{3}^{\ast} - (Y_{2}^{\ast})^{2}) \\ & \hspace{1em} + (((K_{1} - E (X_{1})) Y_{3}^{\ast} - (K_{2} E (X_{2})) Y_{2}^{\ast}) (\mu_{a_{j}} B_{k, 1} + \mu_{b_{j}} A_{k, 1} - \mu_{a_{k}} B_{j, 1} - \mu_{b_{k}} A_{j, 1}) \\ & \hspace{1em} + ((K_{2} - E (X_{2})) Y_{1}^{\ast} - (K_{1} - E (X_{1})) Y_{2}^{\ast}) \\ & \hspace{1em} \times (\mu_{b_{j}} A_{k, 2} + \mu_{a_{k}} B_{j, 2} - \mu_{b_{k}} A_{j, 2} - \mu_{a_{j}} B_{k, 2})) / (Y_{1}^{\ast} Y_{3}^{\ast} - (Y_{2}^{\ast})^{2}) \\ & \hspace{1em} + (((K_{1} - E (X_{1})) Y_{3}^{\ast} - (K_{2} - E (X_{2})) Y_{2}^{\ast}) ((K_{2} - E (X_{2})) Y_{1}^{\ast} - (K_{1} - E (X_{1})) Y_{2}^{\ast}) \\ & \hspace{1em} \times (A_{j, 1} B_{k, 2} + B_{j, 1} A_{k, 2} - A_{j, 2} B_{k, 1} - B_{j, 2} A_{k, 1})) / (Y_{1}^{\ast} Y_{3}^{\ast} - (Y_{2}^{\ast})^{2})^{2} \\ & \hspace{1em} + (((K_{1} - E (X_{1})) Y_{3}^{\ast} - (K_{2} - E (X_{2})) Y_{2}^{\ast})^{2} (B_{j, 1} A_{k, 1} - A_{j, 1} B_{k, 1}) \\ & \hspace{1em} + ((K_{2} - E (X_{2})) Y_{1}^{\ast} - (K_{1} - E (X_{1})) Y_{2}^{\ast})^{2} (A_{j, 2} B_{k, 1} - B_{j, 2} A_{k, 2})) / (Y_{1}^{\ast} Y_{3}^{\ast} - (Y_{2}^{\ast})^{2})^{2}, \\ \\
& E (a_{j}^{2} + b_{j}^{2} \mid \boldsymbol{X} = \boldsymbol{K})
 = (\sigma_{a_{j}}^{2} + \sigma_{b_{j}}^{2}) - ((A_{j, 1}^{2} + B_{j, 1}^{2}) Y_{3}^{\ast} + ((A_{j, 2}^{2} + B_{j, 2}^{2}) Y_{1}^{\ast} \\ & \hspace{1em} - 2 (A_{j, 1} A_{j, 2} - B_{j, 1} B_{j, 2}) Y_{2}^{\ast}) / (Y_{1}^{\ast} Y_{3}^{\ast} - (Y_{2}^{\ast})^{2}) \\ & \hspace{1em} - 2 (((K_{1} - E (X_{1})) Y_{3}^{\ast} - (K_{2} - E (X_{2})) Y_{1}^{\ast}) (\mu_{a{j}} A_{j, 1} - \mu_{b_{j}} B_{j, 1}) \\ & \hspace{1em} + ((K_{2} - E (X_{2})) Y_{1}^{\ast} - (K_{1} - E (X_{1})) Y_{2}^{\ast}) (\mu_{a_{j}} A_{j, 2} + \mu_{b_{j}} B_{j, 2})) / (Y_{1}^{\ast} Y_{3}^{\ast} - (Y_{2}^{\ast})^{2}) \\ & \hspace{1em} + 2 (((K_{1} - E (X_{1})) Y_{3}^{\ast} - (K_{2} - E (X_{2})) Y_{2}^{\ast}) ((K_{2} - E (X_{2})) Y_{1}^{\ast} - (K_{1} - E (X_{1})) Y_{2}^{\ast}) \\ & \hspace{1em} \times (A_{j, 1} A_{j, 2} - B_{j, 1} B_{j, 2})) / (Y_{1}^{\ast} Y_{3}^{\ast} - (Y_{2}^{\ast})^{2})^{2} \\ & \hspace{1em} + (((K_{1} - E (X_{1})) Y_{3}^{\ast} - (K_{2} - E (X_{2})) Y_{2}^{\ast})^{2} (A_{j, 1}^{2} + B_{j, 1}^{2}) \\ & \hspace{1em} + ((K_{2} - E (X_{2})) Y_{1}^{\ast} - (K_{1} - E (X_{1})) Y_{2}^{\ast})^{2} (A_{j, 2}^{2} + B_{j, 2}^{2})) / (Y_{1}^{\ast} Y_{3}^{\ast} - (Y_{2}^{\ast})^{2})^{2}.
\end{align*}
\endgroup
After all the necessary simplifications,
\begin{equation*}
\begin{split}
E (\det (\nabla \boldsymbol{X}) \mid \boldsymbol{X} = \boldsymbol{K})
 & = \sum_{j = 0}^{N} (\sigma_{a_{j}}^{2} + \sigma_{b_{j}}^{2}) \left(\left(\frac{\partial u_{j}}{\partial x}\right)^{2} + \left(\frac{\partial v_{j}}{\partial x}\right)^{2}\right) \\ & \hspace{1em} - \frac{I_{1} S_{1} - I_{2} S_{2} - I_{3} S_{3} - I_{4} S_{4} - I_{5} S_{5}}{Y_{1}^{\ast} Y_{3}^{\ast} - (Y_{2}^{\ast})^{2}},
\end{split}
\end{equation*}
where
\begingroup
\allowdisplaybreaks
\begin{align*}
I_{1}
 & = Y_{3}^{\ast} - \frac{((K_{1} - E (X_{1})) Y_{3}^{\ast} - (K_{2} - E (X_{2})) Y_{2}^{\ast})^{2}}{Y_{1}^{\ast} Y_{3}^{\ast} - (Y_{2}^{\ast})^{2}}, \\
I_{2}
 & = Y_{1}^{\ast} - \frac{((K_{1} - E (X_{1})) Y_{2}^{\ast} - (K_{2} - E (X_{2})) Y_{1}^{\ast})^{2}}{Y_{1}^{\ast} Y_{3}^{\ast} - (Y_{2}^{\ast})^{2}}, \\
I_{3}
 & = Y_{2}^{\ast} + \frac{((K_{1} - E (X_{1})) Y_{3}^{\ast} - (K_{2} - E (X_{2})) Y_{2}^{\ast}) ((K_{2} - E (X_{2})) Y_{1}^{\ast} - (K_{1} - E (X_{1})) Y_{2}^{\ast})}{Y_{1}^{\ast} Y_{3}^{\ast} - (Y_{2}^{\ast})^{2}}, \\
I_{4}
 & = (K_{1} - E (X_{1})) Y_{3}^{\ast} - (K_{2} - E (X_{2})) Y_{2}^{\ast}, \\
I_{5}
 & = (K_{2} - E (X_{2})) Y_{1}^{\ast} - (K_{1} - E (X_{1})) Y_{2}^{\ast}, \\
\end{align*}
\endgroup
and
\begingroup
\allowdisplaybreaks
\begin{align*}
S_{1}
 & = \sum_{j = 0}^{N} \sum_{k = 0}^{N} \left((A_{j, 1} A_{k, 1} + B_{j, 1} B_{k, 1}) \left(\frac{\partial u_{j}}{\partial x} \frac{\partial u_{k}}{\partial x} + \frac{\partial v_{j}}{\partial x} \frac{\partial v_{k}}{\partial x}\right)\right. \\ & \hspace{1em} \left. + (B_{j, 1} A_{k, 1} - A_{j, 1} B_{k, 1}) \left(\frac{\partial v_{j}}{\partial x} \frac{\partial u_{k}}{\partial x} - \frac{\partial u_{j}}{\partial x} \frac{\partial v_{k}}{\partial x}\right)\right) \\
 & = \abs{D_{1}^{\ast}}^{2}, \\
S_{2}
 & = \sum_{j = 0}^{N} \sum_{k = 0}^{N} \left((A_{j, 2} A_{k, 2} + B_{j, 2} B_{k, 2}) \left(\frac{\partial u_{j}}{\partial x} \frac{\partial u_{k}}{\partial x} + \frac{\partial v_{j}}{\partial x} \frac{\partial v_{k}}{\partial x}\right)\right. \\ & \hspace{1em} \left. + (A_{j, 2} B_{k, 2} - B_{j, 2} A_{k, 2}) \left(\frac{\partial v_{j}}{\partial x} \frac{\partial u_{k}}{\partial x} - \frac{\partial u_{j}}{\partial x} \frac{\partial v_{k}}{\partial x}\right)\right) \\
 & = \abs{D_{2}^{\ast}}^{2}, \\
S_{3}
 & = \sum_{j = 0}^{N} \sum_{k = 0}^{N} \left((A_{j, 1} A_{k, 2} + A_{j, 2} A_{k, 1} - B_{j, 1} B_{k, 2} - B_{j, 2} B_{k, 1}) \left(\frac{\partial u_{j}}{\partial x} \frac{\partial u_{k}}{\partial x} + \frac{\partial v_{j}}{\partial x} \frac{\partial v_{k}}{\partial x}\right)\right. \\ & \hspace{1em} \left. + (A_{j, 1} B_{k, 2} + B_{j, 1} A_{k, 2} - A_{j, 2} B_{k, 1} - B_{j, 2} A_{k, 2}) \left(\frac{\partial v_{j}}{\partial x} \frac{\partial u_{k}}{\partial x} - \frac{\partial u_{j}}{\partial x} \frac{\partial v_{k}}{\partial x}\right)\right) \\
 & = \abs{D_{1}^{\ast} + i D_{2}^{\ast}}^{2} - \abs{D_{1}^{\ast}} - \abs{D_{2}^{\ast}}^{2}, \\
S_{4}
 & = \sum_{j = 0}^{N} \sum_{k = 0}^{N} \left((\mu_{a_{j}} A_{k, 1} + \mu_{a_{k}} A_{j, 1} - \mu_{b_{j}} B_{k, 1} - \mu_{b_{k}} B_{j, 1}) \left(\frac{\partial u_{j}}{\partial x} \frac{\partial u_{k}}{\partial x} + \frac{\partial v_{j}}{\partial x} \frac{\partial v_{k}}{\partial x}\right)\right. \\ & \hspace{1em} \left. + (\mu_{b_{k}} A_{j, 2} + \mu_{a_{j}} B_{k, 2} - \mu_{b_{j}} A_{k, 2} - \mu_{a_{k}} B_{j, 2}) \left(\frac{\partial v_{j}}{\partial x} \frac{\partial u_{k}}{\partial x} - \frac{\partial u_{j}}{\partial x} \frac{\partial v_{k}}{\partial x}\right)\right), \\
 & = \abs{M + D_{1}^{\ast}}^{2} - \abs{M}^{2} - \abs{D_{1}^{\ast}}^{2}, \\
S_{5}
 & = \sum_{j = 0}^{N} \sum_{k = 0}^{N} \left((\mu_{a_{j}} A_{k, 2} + \mu_{a_{k}} A_{j, 2} + \mu_{b_{j}} B_{k, 2} + \mu_{b_{k}} B_{j, 2}) \left(\frac{\partial u_{j}}{\partial x} \frac{\partial u_{k}}{\partial x} + \frac{\partial v_{j}}{\partial x} \frac{\partial v_{k}}{\partial x}\right)\right. \\ & \hspace{1em} \left. + (\mu_{b_{k}} A_{j, 2} + \mu_{a_{j}} B_{k, 2} - \mu_{b_{j}} A_{k, 2} - \mu_{a_{k}} B_{j, 2}) \left(\frac{\partial v_{j}}{\partial x} \frac{\partial u_{k}}{\partial x} - \frac{\partial u_{j}}{\partial x} \frac{\partial v_{k}}{\partial x}\right)\right) \\
 & = \abs{M + i D_{2}^{\ast}}^{2} - \abs{M}^{2} - \abs{D_{2}^{\ast}}^{2}.
\end{align*}
\endgroup

After regrouping the terms, with a little algebra we can write
\begin{equation*}
\begin{split}
& E (\det (\nabla \boldsymbol{X}) \mid \boldsymbol{X} = \boldsymbol{K}) 
 = D_{3}^{\ast} - \frac{1}{Y_{1} Y_{3}^{\ast} - (Y_{2}^{\ast})^{2}} (\abs{D_{1}^{\ast}}^{2} I_{1} - \abs{D_{2}^{\ast}}^{2} I_{2} + (\abs{D_{1}^{\ast} + i D_{2}^{\ast}}^{2} - \abs{D_{1}^{\ast}}^{2} - \abs{D_{2}^{\ast}}^{2}) I_{3} \\ & \hspace{2em} - (\abs{M + D_{1}^{\ast}}^{2} - \abs{M}^{2} - \abs{D_{1}^{\ast}}^{2}) I_{4} - (\abs{M + i D_{2}^{\ast}}^{2} - \abs{M}^{2} - \abs{D_{2}^{\ast}}^{2}) I_{5} ).
\end{split}
\end{equation*}
Hence, when combined with the joint density of two normal random variables $X_{1}$ and $X_{2}$ as given by
\begin{equation*}
\begin{split}
& p_{X_{1} X_{2}} (\boldsymbol{K}^{\prime})
 = \frac{1}{2 \pi \sqrt{Y_{1}^{\ast} Y_{3}^{\ast} - (Y_{2}^{\ast})^{2}}} \\ & \hspace{1em}  \times \exp \left(- \frac{(K_{1} - E (X_{1}))^{2} Y_{3}^{\ast} + (K_{2} - E (X_{2}))^{2} Y_{1}^{\ast} - 2 (K_{1} - E (X_{1})) (K_{2} - E (X_{2})) Y_{2}^{\ast}}{2 (Y_{1}^{\ast} Y_{3}^{\ast} - (Y_{2}^{\ast})^{2})}\right),
\end{split}
\end{equation*}
the proof of the theorem is complete.

\section{Ramifications of Theorem \ref{th-4}} \label{sec-5}

Theorem \ref{th-4} has several important consequences. First, if $\mu_{a_{j}} = \mu_{b_{j}} = \mu$ for $0 \leq j \leq N$ in Theorem \ref{th-4}, then the formula for $h_{N, \boldsymbol{K}} (z)$ holds with the following modifications to the auxiliary functions:
\begingroup
\allowdisplaybreaks
\begin{align*}
E (X_{1})
 & = \mu \sum_{j = 0}^{N} (u_{j} - v_{j}), \qquad
E (X_{2})
 = \mu \sum_{j = 0}^{N} (u_{j} + v_{j}), \\
Y_{1}^{\ast} (z)
 & = \sum_{j = 0}^{N} (\sigma_{a_{j}}^{2} u_{j}^{2} + \sigma_{b_{j}}^{2} v_{j}^{2}) - \mu^{2} \left(\sum_{j = 0}^{N} (u_{j} - v_{j})\right)^{2}, \\
Y_{2}^{\ast} (z)
 & = \sum_{j = 0}^{N} (\sigma_{a_{j}}^{2} - \sigma_{b_{j}}^{2}) u_{j} v_{j} - \mu^{2} \left(\sum_{j = 0}^{N} (u_{j} - v_{j})\right) \left(\sum_{j = 0}^{N} (u_{j} + v_{j})\right), \\
Y_{3}^{\ast} (z)
 & = \sum_{j = 0}^{N} (\sigma_{a_{j}}^{2} v_{j}^{2} + \sigma_{b_{j}}^{2} u_{j}^{2}) - \mu^{2} \left(\sum_{j = 0}^{N} (u_{j} + v_{j})\right)^{2},
\end{align*}
\endgroup
and
\begin{equation*}
\begin{split}
A_{j, 1} (z)
 & = \sigma_{a_{j}}^{2} u_{j} - \mu^{2} \sum_{k = 0}^{N} (u_{k} - v_{k}), \\
B_{j, 1} (z)
 & = \sigma_{b_{j}}^{2} v_{j} + \mu^{2} \sum_{k = 0}^{N} (u_{k} - v_{k}),
\end{split}
\qquad
\begin{split}
A_{j, 2} (z)
 & = \sigma_{a_{j}}^{2} v_{j} - \mu^{2} \sum_{k = 0}^{N} (u_{k} + v_{k}), \\
B_{j, 2} (z)
 & = \sigma_{b_{j}}^{2} u_{j} - \mu^{2} \sum_{k = 0}^{N} (u_{k} + v_{k}),
\end{split}
\end{equation*}
for $0 \leq j \leq N$. Observe that $M (z)$ is changed implicitly by the change in $E (\eta_{j})$. If, further, $\mu_{a_{j}} = \mu_{b_{j}} = 0$ for $0 \leq j \leq N$, then $E (X_{1}) = E (X_{2}) = 0$. Thus, $A_{j, 1} = \sigma_{a_{j}}^{2} u_{j}$, $A_{j, 2} = \sigma_{a_{j}}^{2} v_{j}$, $B_{j, 1} = \sigma_{b_{j}}^{2} v_{j}$, $B_{j, 2} = \sigma_{b_{j}}^{2} u_{j}$, and $M (z) = 0$. Hence, Theorem \ref{th-2} is recovered. Further, if $\sigma_{a_{j}}^{2} = \sigma_{b_{j}}^{2} = 1$ for $0 \leq j \leq N$, then Theorem 1 in \cite{CorleyLedoan2020} is recovered.

Second, if $\sigma_{a_{j}}^{2} = \sigma_{b_{j}}^{2} = \sigma^{2}$ for $0 \leq j \leq N$ in Theorem \ref{th-4}, then $Y_{2}^{\ast} (z) = - E (X_{1}) E (X_{2})$ and
\begin{equation*}
Y_{1}^{\ast} (z)
 = \sigma^{2} \sum_{j = 0}^{N} \abs{f_{j} (z)}^{2} - (E (X_{1}))^{2}, \qquad
Y_{3}^{\ast} (z)
 = \sigma^{2} \sum_{j = 0}^{N} \abs{f_{j} (z)}^{2} - (E (X_{2}))^{2}.
\end{equation*}
Observe that
\begin{equation*}
\begin{split}
A_{j, 1} (z) - i B_{j, 1} (z)
 & = \sigma^{2} \overline{f_{j}(z)} - E (X_{1}) E (\eta_{j}), \\
B_{j, 2} (z) - i A_{j, 2} (z)
 & = \sigma^{2} \overline{f_{j}(z)} + i E (X_{2}) E (\eta_{j}).
\end{split}
\end{equation*}
Altogether, the formula for $h_{N, \boldsymbol{K}} (z)$ in Theorem \ref{th-4} now holds with
\begin{equation*}
\begin{split}
D_{0}^{\ast} (z)
 & = \sqrt{\left(\sigma^{2} \sum_{j = 0}^{N} \abs{f_{j} (z)}^{2}\right)^{2} - \sigma^{2} \sum_{j = 0}^{N} \abs{f_{j} (z)}^{2} \left|\sum_{j = 0}^{N} f_{j} (z) E (\eta_{j})\right|^{2}}, \\
D_{1}^{\ast} (z)
 & = \sigma^{2} \sum_{j = 0}^{N} \overline{f_{j}(z)} f_{j}^{\prime} (z) + E (X_{1}) \sum_{j = 0}^{N} f_{j}^{\prime} (z) E (\eta_{j}), \\
D_{2}^{\ast} (z)
 & = \sigma^{2} \sum_{j = 0}^{N} \overline{f_{j}(z)} f_{j}^{\prime} (z) + i E (X_{2}) \sum_{j = 0}^{N} f_{j}^{\prime} (z) E(\eta_{j}), \\
D_{3}^{\ast} (z)
 & = 2 \sigma^{2} \sum_{j = 0}^{N} \abs{f_{j} (z)}^{2}.
\end{split}
\end{equation*}
The form of $M (z)$ remains unchanged. Further, if $\mu_{a_{j}} = \mu_{b_{j}} = \mu$ for $0 \leq j \leq N$,
\begin{equation*}
E (X_{1})
 = \mu \sum_{j = 0}^{N} (u_{j} - v_{j}), \qquad
E (X_{2})
 = \mu \sum_{j = 0}^{N} (u_{j} + v_{j}).
\end{equation*}
Then the formula for $h_{N, \boldsymbol{K}} (z)$ in Theorem \ref{th-4} now holds with
\begingroup
\allowdisplaybreaks
\begin{align*}
D_{0}^{\ast} (z)
 & = \sqrt{\left(\sigma^{2} \sum_{j = 0}^{N} \abs{f_{j} (z)}^{2}\right)^{2} - 2 \mu^{2} \sigma^{2} \sum_{j = 0}^{N} \abs{f_{j} (z)}^{2} \left|\sum_{j = 0}^{N} f_{j} (z)\right|^{2}}, \\
D_{1}^{\ast} (z)
 & = \sigma^{2} \sum_{j = 0}^{N} \overline{f_{j}(z)} f_{j}^{\prime} (z) + \mu \sum_{j = 0}^{N} (u_{j} - v_{j}) \sum_{j = 0}^{N} f_{j}^{\prime} (z) E (\eta_{j}), \\
D_{2}^{\ast} (z)
 & = \sigma^{2} \sum_{j = 0}^{N} \overline{f_{j}(z)} f_{j}^{\prime} (z) + i \mu \sum_{j = 0}^{N} (u_{j} + v_{j}) \sum_{j = 0}^{N} f_{j}^{\prime} (z) E(\eta_{j}), \\
D_{3}^{\ast} (z)
 & = 2 \sigma^{2} \sum_{j = 0}^{N} \abs{f_{j} (z)}^{2}.
\end{align*}
\endgroup

Third, the following also follows from Theorem \ref{th-4}.

\begin{corollary} \label{cor-3}
If the vector $\boldsymbol{K}$ is restricted to a circle of radius $K > 0$, then for all $N > 1$ one has
\begin{equation*}
\begin{split}
& h_{N, \boldsymbol{K}} (z) \\
 & \hspace{1em} = \frac{1}{2 \pi D_{0}^{\ast} (z)} \exp \left(-\frac{(K - E (X_{1}))^{2} Y_{3}^{\ast} (z) + (K - E (X_{2}))^{2} Y_{1}^{\ast} (z) - 2 (K - E (X_{1})) (K - E (X_{2})) Y_{2}^{\ast} (z)}{2 D_{0}^{\ast} (z)^{2}}\right) \\ & \hspace{2em} \times \left\{D_{3}^{\ast} (z) - \frac{\abs{D_{1}^{\ast} (z)}^{2}}{D_{0}^{\ast} (z)^{2}} \left(Y_{3}^{\ast} (z) - \frac{((K - E (X_{1})) Y_{3}^{\ast} (z) - (K - E (X_{2})) Y_{2}^{\ast} (z))^{2}}{D_{0}^{\ast} (z)^{2}}\right) \right. \\ & \hspace{2em} \left. - \frac{\abs{D_{2}^{\ast} (z)}^{2}}{D_{0}^{\ast} (z)^{2}} \left(Y_{1}^{\ast} (z) - \frac{((K - E (X_{1})) Y_{2}^{\ast} (z) - (K - E (X_{2})) Y_{1}^{\ast} (z))^{2}}{D_{0}^{\ast} (z)^{2}}\right) \right. \\ & \hspace{2em} \left. + \left(\frac{\abs{D_{1}^{\ast} (z) + i D_{2}^{\ast} (z)}^{2} - \abs{D_{1}^{\ast} (z)}^{2} - \abs{D_{2}^{\ast} (z)}^{2}}{D_{0}^{\ast} (z)^{2}}\right) \right. \\ &\hspace{2em} \left. \times \left(Y_{2}^{\ast} (z) - \frac{((K - E (X_{1})) Y_{3}^{\ast} (z) - (K - E (X_{2})) Y_{2}^{\ast} (z)) ((K - E (X_{1})) Y_{2}^{\ast} (z) - (K - E (X_{2})) Y_{1}^{\ast} (z))}{D_{0}^{\ast} (z)^{2}}\right) \right. \\ & \hspace{2em} \left. - \left(\frac{\abs{M (z) + D_{1}^{\ast} (z)}^{2} - \abs{M (z)}^{2} - \abs{D_{1}^{\ast} (z)}^{2}}{D_{0}^{\ast} (z)^{2}}\right) ((K - E (X_{1})) Y_{3}^{\ast} (z) - (K - E (X_{2})) Y_{2}^{\ast} (z)) \right. \\ & \hspace{2em} \left. + \left(\frac{\abs{M (z) + i D_{2}^{\ast} (z)}^{2} - \abs{M (z)}^{2} - \abs{D_{2}^{\ast} (z)}^{2}}{D_{0}^{\ast} (z)^{2}}\right) ((K - E (X_{1})) Y_{2}^{\ast} (z) - (K - E (X_{2})) Y_{1}^{\ast} (z))\right\}.
\end{split}
\end{equation*}
\end{corollary}

Then immediate by Corollary \ref{cor-3} is the following result.

\begin{corollary}
If $\boldsymbol{K}$ is the zero vector, then for all $N > 1$ one has
\begingroup
\allowdisplaybreaks
\begin{align*}
h_{N, \boldsymbol{0}} (z)
 & = \frac{1}{2 \pi D_{0}^{\ast} (z)} \exp \left(-\frac{(E (X_{1}))^{2} Y_{3}^{\ast} (z) + (E (X_{2}))^{2} Y_{1}^{\ast} (z) - 2 E (X_{1}) E (X_{2}) Y_{2}^{\ast} (z)}{2 D_{0}^{\ast} (z)^{2}}\right) \\ & \hspace{1em} \times \left\{D_{3}^{\ast} (z) - \frac{\abs{D_{1}^{\ast} (z)}^{2}}{D_{0}^{\ast} (z)^{2}} \left(Y_{3}^{\ast} (z) - \frac{(E (X_{1}) Y_{3}^{\ast} (z) - E (X_{2}) Y_{2}^{\ast} (z))^{2}}{D_{0}^{\ast} (z)^{2}}\right) \right. \\ & \hspace{1em} \left. - \frac{\abs{D_{2}^{\ast} (z)}^{2}}{D_{0}^{\ast} (z)^{2}} \left(Y_{1} - \frac{(E (X_{1}) Y_{2}^{\ast} (z) - E (X_{2}) Y_{1}^{\ast} (z))^{2}}{D_{0}^{\ast} (z)^{2}}\right) + \left(\frac{\abs{D_{1}^{\ast} (z) + i D_{2}^{\ast} (z)}^{2} - \abs{D_{1}^{\ast} (z)}^{2} - \abs{D_{2}^{\ast} (z)}^{2}}{D_{0}^{\ast} (z)^{2}}\right) \right. \\ &\hspace{1em} \left. \times \left(Y_{2}^{\ast} (z) - \frac{(E (X_{1}) Y_{3}^{\ast} (z) - E (X_{2}) Y_{2}^{\ast} (z)) (E (X_{1}) Y_{2}^{\ast} (z) - E (X_{2}) Y_{1}^{\ast} (z))}{D_{0}^{\ast} (z)^{2}}\right) \right. \\ & \hspace{1em} \left. - \left(\frac{\abs{M (z) + D_{1}^{\ast} (z)}^{2} - \abs{M (z)}^{2} - \abs{D_{1}^{\ast} (z)}^{2}}{D_{0}^{\ast} (z)^{2}}\right) (E (X_{1}) Y_{3}^{\ast} (z) - E (X_{2}) Y_{2}^{\ast} (z)) \right. \\ & \hspace{1em} \left. + \left(\frac{\abs{M (z) + i D_{2}^{\ast} (z)}^{2} - \abs{M (z)}^{2} - \abs{D_{2}^{\ast} (z)}^{2}}{D_{0}^{\ast} (z)^{2}}\right) (E (X_{1}) Y_{2}^{\ast} (z) - E (X_{2}) Y_{1}^{\ast} (z))\right\}.
\end{align*}
\endgroup
\end{corollary}

\begin{acknowledgments}
\textup{During the preparation of this work, the first author received support from the Office of the Vice Chancellor for Research and Dean of the Graduate School at The University of Tennessee at Chattanooga. The second author received support from the National Science Foundation under Grant DMS-1852288.}
\end{acknowledgments}

\bibliographystyle{amsplain}

\end{document}